\documentclass[a4paper,11pt]{article}

\usepackage{amsmath, amssymb, amsthm}
\usepackage{graphics}

\usepackage{psfrag, booktabs}
\usepackage{color, epsfig, float}
\usepackage[dvips, bookmarksopen, colorlinks, linkcolor = blue,
urlcolor = red, citecolor = red, menucolor = blue]{hyperref}

\newtheorem{theorem}{Theorem}[section]
\newtheorem{lemma}[theorem]{Lemma}

\newtheorem{propo}[theorem]{Proposition}

\newtheorem{remark}[theorem]{Remark}

\usepackage{setspace}

\newcommand\F{\mathrm{F}}

\newcommand\pr{\Phi_{\rm rel}}

\newcommand\R{\mathbb{R}}

\newcommand\N{\mathbb{N}}

\DeclareMathOperator{\argmin}{arg\, min}
\DeclareMathOperator{\Mo}{\bf M} \DeclareMathOperator{\nr}{\mathcal
N} \DeclareMathOperator{\ra}{\mathcal R}

\newcommand\norm[1]{\|#1\|}
\newcommand\bracket[1]{\langle#1\rangle}
\newcommand\set[1]{\{#1\}}
\newcommand\abs[1]{|#1|}

\begin{document}
%
%
%
%
%

\title{On Steepest-Descent-Kaczmarz methods for regularizing
systems of nonlinear ill-posed equations} \setcounter{footnote}{1}

\author{
A.~De~Cezaro%
\thanks{IMPA, Estr. D. Castorina 110, 22460-320 Rio de Janeiro, Brazil
 \href{mailto:decezaro@impa.br}{\tt decezaro@impa.br}.}
\and M.~Haltmeier
\thanks{Department of Computer Science, University of Innsbruck,
        Technikerstrasse 21a, A-6020 Innsbruck, Austria
        \href{mailto:markus.haltmeier@uibk.ac.at}{\tt \{markus.haltmeier,otmar.scherzer\}@uibk.ac.at}.}
\and
A.~Leit\~ao%
\thanks{Department of Mathematics, Federal University of St. Catarina,
        P.O. Box 476, 88040-900 Florian\'opolis, Brazil
        \href{mailto:aleitao@mtm.ufsc.br}{\tt aleitao@mtm.ufsc.br}.
} \and O.~Scherzer ${}^\ddag$}
\date{\small \today}

\maketitle

\begin{abstract}
We investigate modified steepest descent methods coupled with a
loping Kaczmarz strategy for obtaining stable solutions of nonlinear
systems of ill-posed operator equations. We show that the proposed
method is a convergent regularization method. Numerical tests are
presented for a linear problem related to photoacoustic tomography
and a non-linear problem related to the testing of  semiconductor
devices.
\end{abstract}

\medskip


\section{Introduction} \label{sec:intro}

In this paper we propose a new method for obtaining regularized
approximations of systems of nonlinear ill-posed operator equations.

The \textit{inverse problem} we are interested in consists of
determining an unknown physical quantity $x \in X$ from the set of
data $(y_0, \dots,  y_{N-1}) \in Y^N$, where $X$, $Y$ are Hilbert
spaces and $N \geq 1$. In practical situations, we do not know the
data exactly. Instead, we have only approximate measured data
$y_i^\delta \in Y$ satisfying
\begin{equation}\label{eq:noisy-i}
    \norm{ y_i^\delta - y_i } \le \delta_i \, , \ \ i = 0, \dots, N-1 \, ,
\end{equation}
with $\delta_i > 0$ (noise level). We use the notation $\delta :=
(\delta_0, \dots, \delta_{N-1})$. The finite set of data above is obtained
by indirect measurements of the parameter, this process being described by
the model
\begin{equation}\label{eq:inv-probl}
    \F_{i}(x)  =  y_{i} \, , \ \ i = 0, \dots, N-1 \, ,
\end{equation}
where $\F_i: D_i \subset X \to Y$, and $D_i$ are the corresponding
domains of definition.

Standard methods for the solution of system \eqref{eq:inv-probl} are
based in the use of \textit{Iterative type} regularization methods
\cite{BakKok04, EngHanNeu96, HanNeuSch95, KalNeuSch08, Lan51} or
\textit{Tikhonov type} regularization methods \cite{EngHanNeu96,
Mor93, SeiVog89, Tik63b, TikArs77} after rewriting
(\ref{eq:inv-probl}) as a single equation $\F( x ) =  y$, where
\begin{align} \label{eq:single-op}
    \F  : = ( \F_0, \dots, \F_{N-1} ): \bigcap_{i=0}^{N-1} D_i \to Y^N
\end{align}
and $y := (y_0, \dots,  y_{N-1})$. However these methods become
inefficient if $N$ is large or the evaluations of $\F_i(x)$ and
$\F_i'(x)^\ast$ are expensive. In such a situation, Kaczmarz type
methods \cite{EggHerLen81, Kac37, Mcc77, NatWue01} which cyclically
consider each equation in (\ref{eq:inv-probl}) separately are much
faster \cite{Nat97} and are often the method of choice in practice.

For recent analysis of Kaczmarz type methods for systems of
ill-posed equations, we refer the reader to \cite{BurKal06, HKLS07,
HLS07, KowSch02}.
\medskip

The starting point of our approach is the steepest descent method
\cite{EngHanNeu96, Sch96} for solving ill-posed problems. Motivated
by the ideas in \cite{HKLS07, HLS07}, we propose in this article a
\textit{loping Steepest-Descent-Kaczmarz method} (\textsc{l-SDK}
method) for the solution of (\ref{eq:inv-probl}). This iterative
method is defined by
\begin{equation} \label{eq:lsdk}
    x_{k+1}^\delta = x_{k}^\delta - \omega_k \alpha_k s_k \, ,
\end{equation}
where
\begin{align} \label{eq:def-sk}
s_k &:= \F_{[k]}'(x_{k}^\delta)^*
       ( \F_{[k]}(x_{k}^\delta) - y_{[k]}^\delta ) \, ,
\\[0.5em]
\label{eq:def-omk}
    \omega_k &:=
    \begin{cases}
           1  & \norm{\F_{[k]}(x_{k}^\delta) - y_{[k]}^\delta}
                \geq \tau \delta_{[k]} \\
           0  & \text{otherwise}
         \end{cases} \, ,
\\
\label{eq:def-alk} \alpha_k &:=
     \begin{cases}
      \pr \left( \| s_k \|^2 / \| \F_{[k]}'(x_{k}^\delta) s_k \|^2 \right) & \omega_k = 1 \\
      \alpha_{\rm min}                                                  & \omega_k = 0
   \end{cases} \, .
\end{align}
Here $\alpha_{\rm min} > 0 $, $\tau \in [2,\infty)$ are appropriate
chosen numbers (see \eqref{eq:def-au}, \eqref{eq:def-tau} below),
$[k] := (k \mod N) \in \set{0, \dots, N-1}$, and $x_0^\delta = x_0
\in X$ is an initial guess, possibly incorporating some \textit{a
priori} knowledge about the exact solution. The function $\pr: (0,
\infty) \to (0, \infty)$ defines a sequence of relaxation parameters
and is assumed to be continuous, monotonically increasing, bounded
by a constant $\alpha_{\rm max}$, and to satisfy $\Phi(s) \leq s$
(see Figure \ref{fg:relax}).

\begin{psfrags}
\psfrag{a}{$\alpha$} \psfrag{relax}{$\pr(\alpha)$}
\begin{figure}[tb!]
\centering
\includegraphics[width=0.8\textwidth]{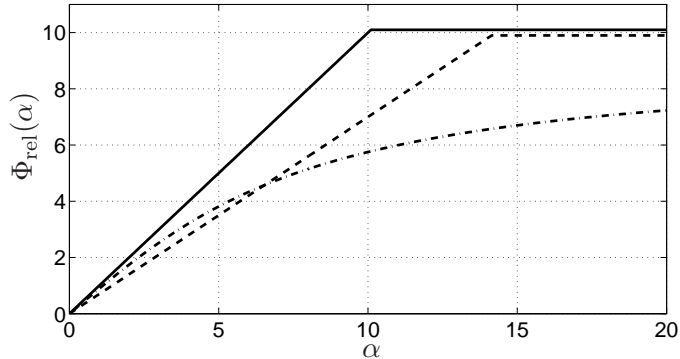}
\caption{Typical examples for relaxation function $\pr$.}
\label{fg:relax}
\end{figure}
\end{psfrags}

If $M$ is an upper bound for $\norm{\F_{[k]}'(x)}$, then $ \| s_k
\|^2 / \| \F_{[k]}'(x_{k}^\delta) s_k \|^2 \geq 1/M^2$ (cf. Lemma
\ref{lem:alk-baf0}). Hence the relaxation function $\pr$ needs only
be defined on $[1/M^2, \infty)$. In particular, if one chooses $\pr(
s ) = \alpha_{\rm min}$ being constant on that interval, then
$\alpha_k = \alpha_{\rm min}$ and the \textsc{l-SDK} method reduces
to the loping Landweber-Kaczmarz (\textsc{l-LK}) method considered
in \cite{HKLS07,HLS07}. The convergence analysis of the
\textsc{l-LK} method requires $\alpha_{\rm min} \leq 1 / M^2$,
whereas the adaptive choice of the relaxation parameters in the
present paper allows $\alpha_k$ being much larger than $1/M^2$.

The \textsc{l-SDK} method consists in incorporating the Kaczmarz
strategy (with the loping parameters $\omega_k$) in the steepest
descent method. This strategy is analog to the one introduced in
\cite{HLS07} regarding the Landweber-Kaczmarz iteration. As usual in
Kaczmarz-type algorithms, a group of $N$ subsequent steps (starting
at some multiple $k$ of $N$) shall be called a {\em cycle}. The
iteration should be terminated when, for the first time, all $x_k$
are equal within a cycle. That is, we stop the iteration at
\begin{equation} \label{eq:def-discrep}
    k_*^\delta  :=  \argmin \set{ l N  \in \N : \, x_{lN}^\delta = x_{lN+1}^\delta = \cdots = x_{lN+N-1}^\delta  } \, ,
\end{equation}
Notice that $k_*^\delta$ is the smallest multiple of $N$ such that
\begin{equation} \label{eq:def-discrep2}
    x_{k_*^\delta} = x_{k_*^\delta+1} = \dots = x_{k_*^\delta+N-1} \, .
\end{equation}
In the case of noise free data, $\delta_i = 0$ in
\eqref{eq:noisy-i}, we choose $\omega_k \equiv 1$ and the iteration
\eqref{eq:lsdk} - \eqref{eq:def-alk} reduces to the {\em
Steepest-Descent-Kaczmarz} (\textsc{SDK}) method, which is closely
related to the {\em Landweber-Kaczmarz} (\textsc{LK}) method
considered in \cite{KowSch02}.

It is worth noticing that, for noisy data, the \textsc{l-SDK} method
is fundamentally different from the \textsc{SDK} method: The
bang-bang relaxation parameter $\omega_k$ effects that the iterates
defined in (\ref{eq:lsdk}) become stationary if \textit{all
components} of the residual vector $\norm{\F_i(x_k^\delta) -
y_i^{\delta}}$ fall below a pre-specified threshold. This
characteristic renders \eqref{eq:lsdk} - (\ref{eq:def-alk}) a
regularization method (see Section~\ref{sec:conv-an}). Another
consequence of using these relaxation parameters is the fact that,
after a large number of iterations, $\omega_k$ will vanish for some
$k$ within each iteration cycle. Therefore, the computational
expensive evaluation of $\F'_{[k]}(x_{k})^*$ might be loped, making
the \textsc{l-SDK} method in \eqref{eq:lsdk} - \eqref{eq:def-alk} a
fast alternative to the \textsc{LK} method in \cite{KowSch02}. Since
in praxis the steepest descent method performs better than the
Landweber method, the \textsc{l-SDK} is expected to be more
efficient than the \textsc{l-LK} method \cite{HKLS07,HLS07}. Our
numerical experiments (mainly for the nonlinear problem considered in
Section~\ref{sec:app-2})  corroborate this conjecture.
\bigskip

The article is  outlined as follows. In Section \ref{sec:basic} we
formulate basic assumptions and derive some auxiliary estimates
required for the analysis. In Section \ref{sec:conv-an} we provide a
convergence analysis for the \textsc{l-SDK} method. In Sections
\ref{sec:app-1} and \ref{sec:app-2} we compare the numerical
performance of the \textsc{l-SDK} method with other standard methods
for inverse problems in photoacoustic tomography and in
semiconductors respectively.

\section{Assumptions and Basic Results} \label{sec:basic}

We begin this section by introducing some assumptions, that are
necessary for the convergence analysis presented in the next
section. These assumptions derive from the classical assumptions
used in the analysis of iterative regularization methods
\cite{EngHanNeu96, KalNeuSch08, Sch96}.

First, we assume that the operators $\F_i$ are continuously
Fr\'echet differentiable, and also that there exist $x_0 \in X$, $M
> 0$, and $\rho > 0$ such that
\begin{equation} \label{eq:a-dfb} \| \F_i'(x) \|
\le M \, , \quad \ x \in B_\rho(x_0) \subset \bigcap_{i=0}^{N-1} D_i
\, .
\end{equation}
Notice that $x_0^\delta =x_0$ is used as starting value of the
\textsc{l-SDK} iteration. Next we make an uniform assumption on the
nonlinearity of the operators $\F_i$. Namely, we assume that the
{\em local tangential cone condition} \cite{EngHanNeu96,
KalNeuSch08}
\begin{equation} \label{eq:a-tcc}
\begin{aligned}
\| \F_i(x) - \F_i(\bar{x}) -  \F_i'(x)( &x - \bar{x}) \|_Y
\\
&\leq
     \eta \norm{ \F_i(x)-\F_i(\bar{x}) }_{Y} \, , \qquad
     x, \bar{x} \in B_{\rho}(x_0)
\end{aligned}
\end{equation}
holds for some $\eta < 1 / 2$. Moreover, we assume the existence of
and element \begin{equation} \label{eq:a-xstar}
    x^* \in B_{\rho/2}(x_0) \,  \text{ such that } \, \F(x^*) = y \,.
\end{equation}
where $y = (y_0, \dots,  y_{N-1})$ are the exact data satisfying
\eqref{eq:noisy-i}.

We are now in position to choose the positive constants $\alpha_{\rm
min}$,  $\tau$ in \eqref{eq:def-alk}, \eqref{eq:def-omk}. For the
rest of this article we shall assume
\begin{align} \label{eq:def-au}
    & \alpha_{\rm min}
    :=  \pr \left( 1/M^2 \right) \,,
    \\ \label{eq:def-tau}
    & \tau
    \,
    \geq 2
    \, \frac{1+\eta}{1-2\eta} \geq 2 \, .
\end{align}
In particular, for linear problems we can choose $\tau$ equal to 2.
\medskip

In the sequel we verify some basic results that are necessary for
the convergence analysis derived in the next section. The first
result concerns the well-definedness and positivity of the
relaxation parameter $\alpha_k$.

\begin{lemma} \label{lem:alk-pos}
Let assumptions \eqref{eq:a-dfb} - \eqref{eq:a-xstar} be satisfied.
Then the coefficients $\alpha_k$ in \eqref{eq:def-alk} are
well-defined and positive.
\end{lemma}

\begin{proof}
If $\omega_k = 0$, the assertion follows from \eqref{eq:def-alk}. If
$\omega_k = 1$, then $\| \F_{[k]}(x_{k}^\delta) - y_{[k]}^\delta\|
\geq \tau \delta_{[k]}$ and the assertion is a consequence of
\cite[Lemma~3.1]{Sch96}, applied to $\F_{[k]}$ instead of $\F$.
\end{proof}

In the next lemma we prove an estimate for the step size of the
\textsc{l-SDK} iteration.

\begin{lemma} \label{lem:sk-estim}
Let $s_k$ and $\alpha_k$  be defined by (\ref{eq:def-sk}) and
(\ref{eq:def-alk}). Then
\begin{equation} \label{eq:sk-estim}
    \alpha_k \norm{s_k}^2
    \ \leq \
    \norm{ \F_{[k]}(x_{k}^\delta) - y_{[k]}^\delta }^2
    \,, \qquad k \in \N  \, .
\end{equation}
\end{lemma}
\begin{proof}
It is enough to consider the case $\omega_k = 1$. It follows from
\eqref{eq:def-alk} that
\begin{equation} \label{eq:aux1}
\alpha_k \norm{s_k}^2 \ =  \ \pr \left( \frac{ \| s_k \|^2}{\|
\F_{[k]}'(x_k^\delta) \, s_k \|^2} \right) \| s_k \|^2
 \ \leq  \ \frac{ \| s_k
\|^4}{\| \F_{[k]}'(x_k^\delta) \, s_k \|^2} \, .
\end{equation}
Moreover, from the definition of $s_k$ we obtain
\begin{align*}
&\| \F_{[k]}'(x_k^\delta) \, s_k \| \ = \ \| \F_{[k]}'(x_{k}^\delta)
\F_{[k]}'(x_{k}^\delta)^*
                   [ \F_{[k]}(x_{k}^\delta) - y_{[k]}^\delta ] \|  \, ,
\\
& \norm{s_k}^2 \ \le \ \| \F_{[k]}'(x_{k}^\delta)
\F_{[k]}'(x_{k}^\delta)^*
                   [ \F_{[k]}(x_{k}^\delta) - y_{[k]}^\delta ] \| \,
                   \| \F_{[k]}(x_{k}^\delta) - y_{[k]}^\delta \|  \, .
\end{align*}
Now, substituting the last two expressions in \eqref{eq:aux1},
shows \eqref{eq:sk-estim}.
\end{proof}

The following Lemma  is an important auxiliary result, which will be
used at several places throughout this article.

\begin{lemma} \label{lem:monot-aux}
Let $x_k^\delta$, $\alpha_k$, $\omega_k$, and $s_k$ be defined by
\eqref{eq:lsdk} - \eqref{eq:def-alk} and assume that
 \eqref{eq:a-dfb} - \eqref{eq:a-xstar} hold true. If $x_k^\delta \in
B_{\rho/2}(x^*)$ for some $k \geq 0$, then
\begin{align} \label{eq:monot-aux}
  \| & x_{k+1}^\delta  - x^* \|^2 -  \| x_k^\delta - x^* \|^2
 \\ \nonumber
 & \le \omega_k \alpha_k \norm{ \F_{[k]}(x_k^\delta) - y_{[k]}^\delta }
  \Big( (2\eta-1) \norm{ \F_{[k]}(x_k^\delta) - y_{[k]}^\delta }
        + 2(1+\eta)\delta_{[k]} \Big) \,. \quad{}
\end{align}
\end{lemma}

\begin{proof}
If $\omega_k = 0$, then $x_{k+1} = x_k$ and \eqref{eq:monot-aux}
follows with equality. If $\omega_k = 1$, it follows from
(\ref{eq:lsdk}) and (\ref{eq:def-sk}) and Lemma~\ref{lem:sk-estim}
that
\begin{align*}
&
\norm{ x_{k+1}^\delta - x^*}^2 - \norm{x_k^\delta - x^*}^2 \\[1ex]
&\quad = 2 \langle x_k^\delta - x^* , \ x_{k+1}^\delta -
x_k^\delta \rangle
  +  \norm{ x_{k+1}^\delta - x_k^\delta}^2
\\
& \quad = 2  \alpha_k
  \langle x_k^\delta - x^* , \
        \F'_{[k]}(x_k^\delta)^* (y_{[k]}^\delta - \F_{[k]}(x_k^\delta) ) \rangle
  + \alpha_k^2 \norm{ s_k }^2
\\
& \quad \leq 2 k \alpha_k
  \langle y_{[k]}^\delta-\F_{[k]}(x_k^\delta),
  \F'_{[k]}(x_k^\delta) (x_k^\delta - x^*) \rangle
  +  \alpha_k \norm{ \F_{[k]}(x_{k}^\delta) - y_{[k]}^\delta }^2 \\
  & \quad
\leq \alpha_k \Bigl(
  2 \langle y_{[k]}^\delta  -  \F_{[k]}(x_k^\delta), \
          \F'_{[k]}(x_k^\delta)
          (x_k^\delta - x^*) - \F_{[k]}(x^*) + \F_{[k]}(x_k^\delta)  \rangle \\
& \qquad\qquad\qquad  + 2
  \langle y_{[k]}^\delta  -  \F_{[k]}(x_k^\delta), \
          y_{[k]} - y_{[k]}^\delta \rangle
         -
         \norm{  y_{[k]}^\delta -  \F_{[k]}(x_k^\delta)}^2 \; \Bigr)
         \,.
\end{align*}
Now, applying \eqref{eq:a-tcc} with $x = x^*$ and $\bar x
=x_k^\delta \in B_{\rho/2}(x^*) \subset B_\rho(x_0)$, leads to
\begin{align*}
&
\norm{ x_{k+1}^\delta - x^*}^2 - \norm{x_k^\delta - x^*}^2 \\
& \quad \le \omega_k \alpha_k
  \norm{ \F_{[k]}(x_k^\delta) - y_{[k]}^\delta }\,
  \Big( 2\eta \norm{ \F_{[k]}(x_k^\delta) - y_{[k]} }
        + 2\delta_{[k]}
        - \norm{\F_{[k]}(x_k^\delta) - y_{[k]}^\delta} \, \Big) \,.
\end{align*}
The last inequality and \eqref{eq:noisy-i} show \eqref{eq:monot-aux}.%
\end{proof}

Our next goal is to prove a monotony property, known to be satisfied
by other iterative regularization methods, e.g., by the Landweber
\cite{EngHanNeu96}, the steepest descent \cite{Sch96}, the
\textsc{LK} \cite{KowSch02}, and the \textsc{l-LK} \cite{HLS07}
method.

\begin{propo}[Monotonicity]\label{pr:mon}
Under the assumptions of Lemma~\ref{lem:monot-aux},
\begin{align}\label{eq:sdk-monot1}
\| x_{k+1}^\delta - x^* \|^2 & \le  \| x_k^\delta - x^* \|^2 \, , \
   \qquad  k \in  \N \,.
\end{align}
Moreover, all iterates $x_k^\delta$ remain in $B_{\rho/2}(x^*)
\subset B_{\rho}(x_0)$ and satisfy \eqref{eq:monot-aux}.
\end{propo}

\begin{proof}
From \eqref{eq:a-xstar} it follows that $x_0 \in B_{\rho/2}(x^*)$.
If $\omega^\delta = 0$, then $x_1$ satisfies \eqref{eq:sdk-monot1}
with equality and $x_1 \in B_{\rho/2}(x^*) \subset B_{\rho}(x_0)$.
If $\omega^\delta \not= 0$, then Lemma~\ref{lem:monot-aux} implies
\begin{align*}
\norm{ x_{1}^\delta - x^*}^2 - \norm{x_0^\delta - x^*}^2 & \geq
    (2\eta-1) \norm{ \F_{0}(x_{0}^\delta) - y^{\delta,0} }
          + 2(1+\eta)\delta^{0}
\\
& \geq  \delta_0 \Big( (2\eta-1) \tau + 2(1+\eta) \Big) \geq 0 \,.
\end{align*}
Therefore \eqref{eq:sdk-monot1}, for $k=0$, follows from
\eqref{eq:def-tau}. In particular, $x_1 \in B_{\rho/2}(x^*) $. An
inductive argument implies \eqref{eq:sdk-monot1} and that $x_{k} \in
B_{\rho/2}(x^*) \subset B_{\rho}(x_0)$ for all $k \in \N$. The
assertions therefore follows from Lemma~\ref{lem:monot-aux}.
\end{proof}

\section{Convergence Analysis of the Loping Steepest Descent Kaczmarz
Method} \label{sec:conv-an}

In this section we provide a complete convergence analysis for the
\textsc{l-SDK} iteration, showing that it is a convergent
regularization method in the sense of \cite{EngHanNeu96} (see
Theorems~\ref{th:exact} and~\ref{th:noisy} below). Throughout this
section, we assume that (\ref{eq:a-dfb}) - (\ref{eq:def-tau}) hold,
and that $x_k^\delta$, $\alpha_k$, $\omega_k$, and $s_k$ are defined
by \eqref{eq:lsdk} - \eqref{eq:def-alk}.
\medskip

Our first goal is to prove convergence of the \textsc{l-SDK}
iteration for $\delta = 0$. For exact data $y = (y_0, \dots,
y_{N-1})$, the iterates in \eqref{eq:lsdk} are denoted by
$x_k$.\footnote{This is a standard notation used in the literature.}

\begin{lemma} \label{lem:min-norm}
There exists an $x_0$-minimal norm solution of \eqref{eq:inv-probl} in
$B_{\rho/2}(x_0)$, i.e., a solution $x^\dag$ of \eqref{eq:inv-probl} such that
\begin{align*}
    \norm{ x^\dag - x_0 }
    =
    \inf
    \big\{ \norm{x - x_0}  : x \in B_{\rho/2}(x_0) \text{ and } \F(x) = y \}
    \,.
\end{align*}
Moreover, $x^\dag $ is the only solution of \eqref{eq:inv-probl} in
$B_{\rho/2}(x_0) \cap \bigl( x_0  +  \ker(\F '( x^\dag ))^\perp
\bigr) $.
\end{lemma}

\begin{proof}
Lemma \ref{lem:min-norm} is a consequence of
\cite[Proposition~2.1]{HanNeuSch95}. A detailed proof can be found
in \cite{KalNeuSch08}.
\end{proof}

\begin{lemma} \label{lem:alk-baf0}
For all $k \in \N$, we have $\alpha_k \geq \alpha_{ \rm min}$.
\end{lemma}

\begin{proof}
For $\omega_k = 0$ the claimed estimate holds with equality. If
$\omega_k = 1$, it follows from \eqref{eq:a-dfb} that
\begin{align*}
    \norm{ s_k }^2 / \norm{ \F_{[k]}'(x_k^\delta) \, s_k }^2 \ \ge \
    \norm{ \F_{[k]}'(x_k^\delta) }^{-2}
    \ \geq \  1/M^2 \,.
\end{align*}
Now the monotonicity of $\pr$ implies $\alpha_k \ge \pr(M^{-2})
=\alpha_{\rm min}$.
\end{proof}

Throughout the rest of this article,  $x^\dag$ denotes the
$x_0$-minimal norm solution of \eqref{eq:inv-probl}. We define $e_k
:= x^\dag - x_k$. From Proposition \ref{pr:mon} it follows that
\eqref{eq:monot-aux} holds for all $k$. By summing over all $k$,
this leads to
\begin{align}\label{eq:sum-1}
    \sum_{i=0}^{\infty}
    \alpha_i
    \norm{ y_{[i]} - \F_{[i]}(x_i) }^2
    \ \leq \
    \frac{\norm{x_0 - x^\dag}^2}{1-2\eta}
    \ < \ \infty \,.
\end{align}
Equation \eqref{eq:sum-1} and the monotony of $\norm {e_k}$ shown in
Proposition \ref{pr:mon} are main ingredients in the following proof
of the convergence of the \textsc{SDK} iteration.

\begin{theorem}[Convergence for Exact Data] \label{th:exact}
For exact data, the iteration $(x_k)$ converges to a solution of
\eqref{eq:inv-probl}, as $k \to \infty$. Moreover, if
\begin{equation} \label{eq:kern-cond}
    \nr( \F'(x^\dag) ) \subseteq \nr( \F' (x) )
    \quad \text{ for all } x \in B_\rho(x_0) \,,
\end{equation}
then $x_k \to x^\dag$.
\end{theorem}

\begin{proof}
From \eqref{eq:sdk-monot1} it follows  that $\norm{e_k}$ decreases
monotonically and therefore that $\norm{e_k}$ converges to some
$\epsilon \geq 0$. In the following we show that $e_k$ is in fact a
Cauchy sequence.

For $k = k_0 N + k_1$ and $l = l_0 N + l_1$ with $k \leq l$ and
$ k_1, l_1 \in \set{0, \dots, N-1}$, let $n_0 \in \set{k_0, \dots, l_0}$
be such that
\begin{equation} \label{eq:l1}
    \sum_{i_1 = 0}^{N-1}  \norm{ F_{i_1} (x_{N n_0 + i_1}) - y_{i_1} }
    \leq
    \sum_{i_1=0}^{N-1}  \norm{ F_{i_1} (x_{N i_0 + i_1}) - y_{i_1} }
    \, , \quad i_0 \in \set{k_0, \dots, l_0} \, .
\end{equation}
Then, with $n := N n_0 + N-1$, we  have
\begin{equation}
    \norm{ e_k - e_l }
    \leq
    \norm{ e_k - e_n } +
    \norm{ e_l - e_n }
\end{equation}
and
\begin{equation} \label{eq:conv-h}
\begin{aligned}
    \norm{ e_n - e_k }^2
    & = \norm{ e_k}^2 - \norm{ e_n }^2 + 2\bracket{ e_n - e_k, e_n } \,.
    \\
    \norm{ e_n - e_l }^2
    & = \norm{ e_l}^2 - \norm{ e_n }^2 + 2\bracket{ e_n - e_l, e_n } \,,
\end{aligned}
\end{equation}
For $k \to \infty$, the first two terms of \eqref{eq:conv-h}
converge to $\epsilon - \epsilon =0$. Therefore, in order to show
that $e_k$ is a Cauchy sequence, it is sufficient to prove that
$\langle e_n - e_k , e_n \rangle$ and $\langle e_n - e_l , e_n
\rangle$ converge to zero as $k \to \infty$.

To that end, we write  $i = N i_0 + i_1$, $i_1\in \set{ 0, \dots,
N-1}$ and set $i^* := N n_0 + i_1$. Then, using the definition of
the steepest descent Kaczmarz iteration it follows that
\begin{align}\nonumber
    | \langle & e_n - e_k , e_n \rangle |
    \\ \nonumber
    &
    =
    \biggl|
        \sum_{i = k}^{n-1}
        \alpha_i
        \big\langle
            \F_{i_1}'(x_{i}) ^* \big(y_{i_1} - \F_{i_1}(x_i)\big),  x^\dag  -  x_{n}
        \big\rangle
    \biggr|
    \\ \nonumber
    &\leq \sum_{i= k}^{n-1}
        \alpha_i
        \left|
            \big\langle y_{i_1} - \F_{i_1}(x_i),
                \F_{i_1}'(x_{i}) ( x^\dag  - x_{i^*} ) + \F_{i_1}'(x_{i}) (x_{i^*}-x_n)
            \big\rangle
        \right|
   \\ \nonumber
    & \leq \sum_{i= k}^{n-1}
        \alpha_i  \| y_{i_1} - \F_{i_1}(x_i)\|\,   \|\F_{i_1}'(x_{i}) (x^\dag -  x_{i^*}) \|
   \\[-1.5ex] \label{eq:est-f0}
   & \qquad \qquad \; +\sum_{i= k}^{n-1}
        \alpha_i  \| y_{i_1} - \F_{i_1}(x_i)\|\,   \|\F_{i_1}'(x_{i}) (x_{i^*} -  x_{n}) \|
\end{align}
From \eqref{eq:a-tcc} it follows  that
\begin{align}\label{eq:est-f1}
    \|\F_{i_1}'(x_{i}) (x^\dag -  x_{i^*}) \|
    \leq 2(1+\eta)\| y_{i_1} - \F_{i_1}(x_{i})\|
    + (1+\eta)\| y_{i_1} - \F_{i_1}(x_{i^*})\|\,.
\end{align}
Again using the definition of the steepest descent Kaczmarz
iteration and equations \eqref{eq:def-alk}, \eqref{eq:a-dfb}, it
follows that
\begin{align} \nonumber
      \|  \F_{i_1}'(x_{i}) (x_{i^*} -  &x_{n}) \|
      \leq
      M \, \| x_{i^*} -  x_{n} \|
    \\[0.5ex] \nonumber
    &
    \leq
    M \,
    \sum_{j = i_1}^{N-2}
        \alpha_j \norm{ \F_{j}'(x_{Nn_0+j})^* \big( \F_{j}(x_{Nn_0+j}) - y_j \big)  }
    \\ \label{eq:est-f2}
    &
    \leq
    \alpha_{\rm max} M^2
    \sum_{j = 0}^{N-1}
    \norm{  F_{j}(x_{Nn_0+j}) - y_j }\,.
\end{align}
Substituting \eqref{eq:est-f1}, \eqref{eq:est-f2} in
\eqref{eq:est-f0} leads to
\begin{align*}
    | \langle & e_n - e_k , e_n \rangle |
    \\&
    \leq
    c
    \sum_{i_0 = k_0}^{n-1}
    \sum_{i_1 =0}^{N-1}
    \norm{ y_{i_1} - \F_{i_1}(x_{N i_0 + i_1})}
    \left( \sum_{j = 0}^{N-1}\norm{\F_{j}(x_{N n_0+j}) - y_j } \right)
    \\
    &
    \leq
    c
    \sum_{i_0 = k_0}^{n-1}
    \Bigl( \sum_{i_1 =0}^{N-1} \norm{ y_{i_1} - \F_{i_1}(x_{N i_0 + i_1})}
    \Bigr)^2
    \end{align*}
with  $c := \alpha_{\rm max}(3 + 3\eta + \alpha_{\rm max} M^2)$.
Here we made use of (\ref{eq:l1}). So, we finally obtain the
estimate
\begin{align*}
    | \langle & e_n - e_k , e_n \rangle |
    \leq
    \frac{N c}{\alpha_{\rm min}}
    \sum_{i_0 = k_0}^{n-1}
    \sum_{i_1 = 0}^{N-1}
    \alpha_{N i_0 + i_1} \norm{ y_{i_1} - \F_{i_1}(x_{N i_0 + i_1})}^2 \,.
\end{align*}
Because of \eqref{eq:sum-1}, the last sum tends to zero for $k= N
k_0 + k_1 \to \infty$, and therefore $\bracket{e_n, e_n -e_k} \to
0$.  Analogously one shows that $\bracket{e_n, e_n -e_l} \to 0$.
Therefore $e_k$ is a Cauchy sequence and  $x_k =  x^\dag - e_k$
converges to an element $x^\ast \in X$. Because all residuals
$\norm{\F_{[k]} (x_k) - y_{[k]}}$  tend to zero, $x^*$ is solution
of  \eqref{eq:inv-probl}.
\medskip

Now assume $\nr( \F'(x^\dag) ) \subseteq \nr( \F(x) )$, for $x \in
B_\rho(x_0)$. Then from the definition of $x_k$ it follows that
\begin{align*}
x_{k+1} - x_k  \in  \ra(\F_{[k]}'(x_{k})^*) \subset
    \nr( \F_{[k]}'(x_{k}) )^\bot \subset
    \nr( \F'(x_{k}) )^\bot \subset \nr( \F'(x^\dag) )^\bot
    \,.
\end{align*}
An inductive argument shows that all iterates $x_{k}$ are elements
of $x_0 + \nr( \F'(x^\dag) )^\bot$. Together with the  continuity of
$\F'(x^\dag)$ this implies that $x^* \in x_0 + \nr( \F'(x^\dag)
)^\bot$. By Lemma \ref{lem:min-norm}, $x^\dag$ is the only solution
of \eqref{eq:inv-probl} in $B_{\rho/2}(x_0) \cap \bigl( x_0  +
\nr(\F '( x^\dag ))^\perp \bigr)$, and so the  second assertion
follows.
\end{proof}
\medskip

The second goal in this section is to prove that $x_{k*}^\delta$
converges to a solution of \eqref{eq:inv-probl}, as $\delta \to 0$.
First we verify that, for noisy data, the stopping index
$k_*^\delta$ defined in \eqref{eq:def-discrep} is finite.

\begin{propo}[Stopping Index] \label{prop:st-f}
Assume $\delta_{\rm min} := \min \{ \delta_0, \dots , \delta_{N-1} \}
> 0$. Then $k_*^\delta$ defined in \eqref{eq:def-discrep} is finite,
and
\begin{equation} \label{eq:sdk-monot-res}
\norm{ \F_{i}(x_{k_*^\delta}^\delta) - y_i^\delta} < \tau \delta_i
\,, \qquad i = 0, \dots, N-1 \, .
\end{equation}
\end{propo}

\begin{proof}
Assume that for every $l  \in N$, there exists $i(l) \in
\set{0,\dots, N-1}$ such that $x_{lN + i(l)} \not= x_{lN}$. From
Proposition~\ref{pr:mon} follows that we can apply
\eqref{eq:monot-aux} recursively for $k=1, \dots, lN$ and obtain
\begin{align*} \label{eq:st-finite-aux1}
    - \norm{x_0 - x^*}^2
    \leq  & \sum_{k=1}^{l N} \omega_k \alpha_k
    \norm{ \F_{[k]}(x_{k}^\delta) - y_{[k]}^\delta }
    \\
    & \left(
    2(1+\eta)\delta_{[k]} - (1-2\eta)
    \norm{ \F_{[k]}(x_{k}^\delta) - y_{[k]}^\delta } \right) \,, \quad l \in \N\,.
    \end{align*}
Using the fact that either $\omega_k=0$ or $\norm{
\F_{[k]}(x_{k}^\delta) - y_{[k]}^\delta } \geq \tau \delta_{[k]}$,
we obtain
\begin{equation} \label{eq:st-finite-aux1}
    \norm{x_0 - x^*}^2
    \geq
    \Bigl( \tau (1-2\eta ) - 2(1+\eta) \Bigr)
    \sum_{k=1}^{l N} \,  \omega_{k}  \alpha_k \delta_{[k]}
    \norm{ \F_{[k]}(x_{k}^\delta) - y_{[k]}^\delta }   \, .
\end{equation}
Equation  (\ref{eq:st-finite-aux1}), Lemma \ref{lem:alk-baf0} and
the fact that $x_{l'N + i(l')} \not= x_{l'N}$ for all $l' \in \N$,
imply
\begin{equation} \label{eq:st-finite-aux2}
    \norm{x_0 - x^*}^2
    \ \geq \
    \Bigl( \tau (1-2\eta ) - 2(1+\eta) \Bigr)
    \, l \, \alpha_{\rm min} \, \delta_{\rm min} \, (\tau \delta_{\rm min})
    \, , \quad l \in \N \, .
\end{equation}
The right hand side of (\ref{eq:st-finite-aux2}) tends to infinity,
which gives a contradiction. Consequently, $\{ l  \in \N: x_{lN + i}
= x_{lN} \,, \, 0 \leq i \leq N-1 \}\not =\emptyset$ and the infimum
in (\ref{eq:def-discrep}) takes a finite value.
\medskip

To prove (\ref{eq:sdk-monot-res}), assume to the contrary, that
$\norm{ \F_{i}(x_{k_*^\delta}^\delta) - y_i^\delta} \geq \tau
\delta_i$ for some $i \in \{0, \dots, N-1 \}$. From
\eqref{eq:def-omk} and \eqref{eq:def-discrep} it follows that,
$\omega_{k_*^\delta} = 1$ and $x_{k_*^\delta + i}^\delta = x_{k_*^\delta + i
+ 1}^\delta$ respectively. Thus, Proposition~\ref{pr:mon} and
Lemma~\ref{lem:alk-pos} imply
\begin{align*}
0
    \le (2\eta-1) \norm{ \F_{i}(x_{k_*^\delta}^\delta) - y_i^\delta }
          + 2(1+\eta)\delta_i
    <  \delta_i \Big( (2\eta-1) \tau + 2(1+\eta) \Big) \,.
\end{align*}
This contradicts \eqref{eq:def-tau}, concluding the proof of
\eqref{eq:sdk-monot-res}.
\end{proof}

The last auxiliary result concerns the continuity of $x_{k}^\delta$
at $\delta =0$. For $y, y^\delta  \in Y^N$, $\delta >0$, and $k \in
\N$ we define
\begin{align*}
    \Delta_k(\delta, y, y^\delta)
    &
    :=
    \omega_k \F_{[k]}'(x_k^\delta)^*\bigl( \F_{[k]}(x_k^\delta) - y_{[k]}^\delta\bigr)
    -
    \F_{[k]}'(x_k)^*\bigl( \F_{[k]}(x_k) - y_{[k]} \bigr) \,.
\end{align*}

\begin{lemma} \label{le:cont}
For all $k \in \N$,
\begin{align}\label{eq:cont}
\lim_{\delta \to 0}
     \sup
    \left\{ \norm{ \Delta_k(\delta, y, y^\delta) } :
    y^\delta  \in Y^N, \,\norm{y_{i}- y_{i}^\delta} \leq \delta_{i} \right\}
    = 0
\;.
\end{align}
Moreover,  $x_{k+1}^\delta \to  x_{k+1}$, as $\delta \to 0$.
\end{lemma}

\begin{proof}

We prove Lemma \ref{le:cont} by induction. The case $k=0$ is similar
to the general case and is omitted.

Now, assume $k > 0$ and that (\ref{eq:cont}) holds for all $k' < k$.
First we note that \eqref{eq:cont} and the continuity of $\pr$
obviously imply $x_{k+1}^\delta \to  x_{k+1}$, as $\delta \to 0$.
For the proof of  \eqref{eq:cont}  we consider two cases.
In the first case,  $\omega_k =1$, we have
 \begin{align*}
        \norm{ \Delta_k(\delta, y, y^\delta) }
        &
        = \norm{ \F_{[k]}'(x_k^\delta)^*\bigl( \F_{[k]}(x_k^\delta) - y_{[k]}^\delta\bigr)
        - \F_{[k]}'(x_k)^*\bigl( \F_{[k]}(x_k) - y_{[k]} \bigr)} \,.
\end{align*}
In the second case, $\omega_k =0$, we have
$\norm{\F_{[k]}(x_k^\delta) - y_{[k]}^\delta} \leq \tau \delta^k$
and consequently
    \begin{align*}
    \|\Delta_k&(\delta, y, y^\delta)\|
    \leq
    \norm{\F_{[k]}'(x_k)^* ( \F_{[k]}(x_k) - y_{[k]} )}
     \\
    &\leq
    \norm{\F_{[k]}'(x_k^\delta)}
    \left(
    \norm{ \F_{[k]}(x_k) - \F_{[k]}(x_k^\delta)}
    +
    \norm{ \F_{[k]}(x_k^\delta) - y_k^\delta}
    +
    \norm{y_k^\delta - y_{[k]}}
    \right)
    \\
    & \leq
    \norm{\F_{[k]}'(x_k^\delta)}
    \left(
    \norm{ \F_{[k]}(x_k) - \F_{[k]}(x_k^\delta)}
    +
    (\tau + 1 ) \, \delta_{[k]}
    \right)\,.
    \end{align*}
Now \eqref{eq:cont} follows from \eqref{eq:a-dfb}, the continuity of
$F_{[k]}$ and $F'_{[k]}$, and the induction hypothesis (which
implies $x_k^\delta \to x_k$).
\end{proof}

\begin{theorem}[Convergence for Noisy Data] \label{th:noisy}
Assume  $( \delta^j_{0}, \dots, \delta^j_{N-1} ) $ is a sequence in
$(0,\infty)^N$ with $\lim_{j\to\infty} \delta^j_i = 0$. Let
$(y^j_{0}, \dots, y^j_{N-1})$ be a sequence of noisy data satisfying
\begin{align*}
    \norm{y_i^{j} - y_i} \le \delta^j_i
    \,,
    \qquad
    i=0,\dots,N- 1 \,,\, j \in \N\,,
\end{align*}
and let $k^j := k_*(\delta^j, y^j) $  denote the corresponding
stopping index defined in \eqref{eq:def-discrep}. Then
$x_{k^j}^{\delta^j}$ converges to a solution of
\eqref{eq:inv-probl}, as $j \to \infty$. Moreover, if
\eqref{eq:kern-cond} holds, then $x_{k^j}^{\delta^j} \to x^\dag$.
\end{theorem}

\begin{proof}
Let $x^*$ denote the limit of the iterates $x_k$ which is a solution
of  \eqref{eq:inv-probl}, cf.  Theorem \ref{th:exact}. From  Lemma
\ref{le:cont} and the continuity of $\F_i$ we know that, for any
fixed $k \in \N$,
\begin{align} \label{eq:est-conv}
    x^{\delta^j}_{k} \to x_{k} \,,
    \quad
    \F_{i} ( x_k^{\delta^j} ) \to  \F_{i}(x_{k})
    \,, \qquad \text{ as } j \to \infty \,.
\end{align}
To show that $x_{k^j}^{\delta^j} \to x^*$, we first assume that
$k^j$ has a finite accumulation point $k_*$. Without loss of
generality  we may assume that $k^j = k_*$ for all $j \in \N$.
From Proposition \ref{prop:st-f} we know that $\norm{y_i^{\delta^j} -
\F_{i}(x^{\delta^j}_{k_*})} < \tau \delta_{i}^{j}$ and, by taking the
limit $j \to \infty$, that $ y_{i} = \F_{i}(x_{k_*}) $. Consequently
$x_{k_*} = x^*$ and $x^j_{k_*} \to x^{*}$ as
 $j \to \infty$.

It remains to consider the case where $k^j \to \infty$ as $j \to
\infty$. To that end let $\varepsilon >0$. Without loss of
generality we assume that $k^j$ is monotonically increasing.
According to Theorem \ref{th:exact}  we can choose $n \in \N$ such
that $\norm{x_{k^n}  -  x^*} < \varepsilon/2$. Equation
\eqref{eq:est-conv} implies that there exists $j_0  >  n$ such that
$\norm{x^{\delta^{j}}_{k^n} - x_{k^n}} < \varepsilon/2$ for all $j
\geq j_0$.  This and Proposition \ref{pr:mon} imply
\begin{align*}
    \norm{x^{\delta^{j}}_{k^j}&  - x^*}
    \leq
    \norm{x^{\delta^{j}}_{k^n}  - x^*}
    \\
    & \leq
    \norm{x^{\delta^{j}}_{k^n}  - x_{k^n}}
    +
    \norm{x_{k^n}  - x^*}
    < \frac{\varepsilon}{2} + \frac{\varepsilon}{2} = \varepsilon
    \,, \quad \text{ for } j \geq j_0 \,.
\end{align*}
Consequently, $x^{\delta^{j}}_{k^j} \to x^*$.
\medskip

If \eqref{eq:kern-cond} holds true,  then by  Theorem
\ref{th:exact}, $x^* =  x^\dag$. Therefore $x^{\delta^{j}}_{k^j} \to
x^\dag$, which concludes the proof.
\end{proof}

\begin{remark}
In  standard iterative regularization methods the number of
performed iterations plays the role of the regularization parameter
\cite{EngHanNeu96, KalNeuSch08}. A parameter choice rule corresponds
to the choice of an appropriate stopping index $k_*^\delta =
k(\delta, y^\delta)$.

For the loping Kaczmarz iterations analyzed in this article, the
situation is quite different. If $k$ is fixed, then the iterates
$x_k^\delta$, do \emph{not}  depend continuously on data
$y_i^\delta$. However, for a fixed sequence $(\omega_k)$ of loping
parameters, the iterates $x_k^\delta$ \emph{do depend} continuously
on $y_i^\delta$: Now, the loping sequences $(\omega_k)$ play the
role of the regularization parameters and  the particular sequence
$\omega_k = \omega_k(\delta, y^\delta)$, depending on $\delta_i$ and
the noisy data $y_i^\delta$,  is the  a-posteriori \emph{parameter
choice} rule.
\end{remark}

\section{Limited View Problem in Photoacoustic Computed Tomography}
\label{sec:app-1}

In this section we compare the numerical performance of loping
Kaczmarz methods  applied to a system of linear equations related to
a limited view problem in {\em photoacoustic computed tomography}
\cite{FinRak07, KucKun07, PalNusHalBur07, XuMWan06}.
\medskip

Let $X := L^2(D)$ denote the Hilbert space of all square integrable
functions in the unit disc $D \subset \R^2$, and let $Y$ denote the
Hilbert space of all functions $y: [0,2] \to \R$ with $ \norm{y}^2
:= \int_0^2 y(t) t dt <  \infty $. We consider the system
\begin{equation}\label{eq:ta:eq}
    \Mo_i x
    =
    y_i
    \,,
    \quad \;  i = 0, \dots, N-1\;,
\end{equation}
where  $\Mo_i: X \to Y$,
\begin{equation} \label{eq:smean}
\begin{aligned}
    (\Mo_i x)(t)
   : =
    \frac{1}{\sqrt{\pi}}
        \int_{S^1}
        x(\xi_{i} + t \sigma)
        \, d \Omega(\sigma)\,, \qquad t \in [0,2]\,,
\end{aligned}
\end{equation}
correspond to a scaled version of the \textit{circular mean Radon
transform}. Solving (\ref{eq:ta:eq}) is the crucial step in
three-dimensional photoacoustic computed tomography with integrating
linear detectors \cite{HalFid07,PalNusHalBur07}, where the centers
of integration, $\xi_{i}$, correspond to the positions of the linear
detectors. We are particularly interested in the incomplete data
case (limited view problem), where the centers $\xi_i = \bigl(
\sin(\pi i /(N-1)) \cos(\pi i/(N-1)) \bigr)$ are uniformly
distributed on the semicircle $S^1_+ := \set{ \xi = (\xi^1, \xi^2)
\in \partial D:  \xi^1 \geq 0}$. Micro-local analysis predicts, that
if the  centers do not cover the whole circle, certain details (the
\emph{invisible} boundaries) of $x$ outside the detection region
(convex hull of $S^1_+$) cannot be recovered \cite{LoiQui00,
PalNusHalBur07b, XuYEtAl04}.

The operators $\Mo_i$ are linear, bounded, and satisfy $\norm{\Mo_i}
\leq 1$, \cite{HKLS07}. For linear bounded operators, the
\emph{tangential cone condition}~(\ref{eq:a-tcc}) is  satisfied with
$\eta =0$. Consequently, the analysis of Section \ref{sec:conv-an}
applies, and  the \textsc{l-SDK} method \eqref{eq:lsdk} -
\eqref{eq:def-alk} provides  a convergent regularization method for
solving (\ref{eq:ta:eq}). The adjoint of $\Mo_i$, required in
\eqref{eq:def-sk} is given by $(\Mo_i^\ast y) ( \xi )  = y(\abs{
\xi_i - \xi})/ \sqrt \pi$, \cite{HKLS07}.

\begin{figure}[t!]
\centering
\includegraphics[height = 0.4\textwidth, width=0.45\textwidth]{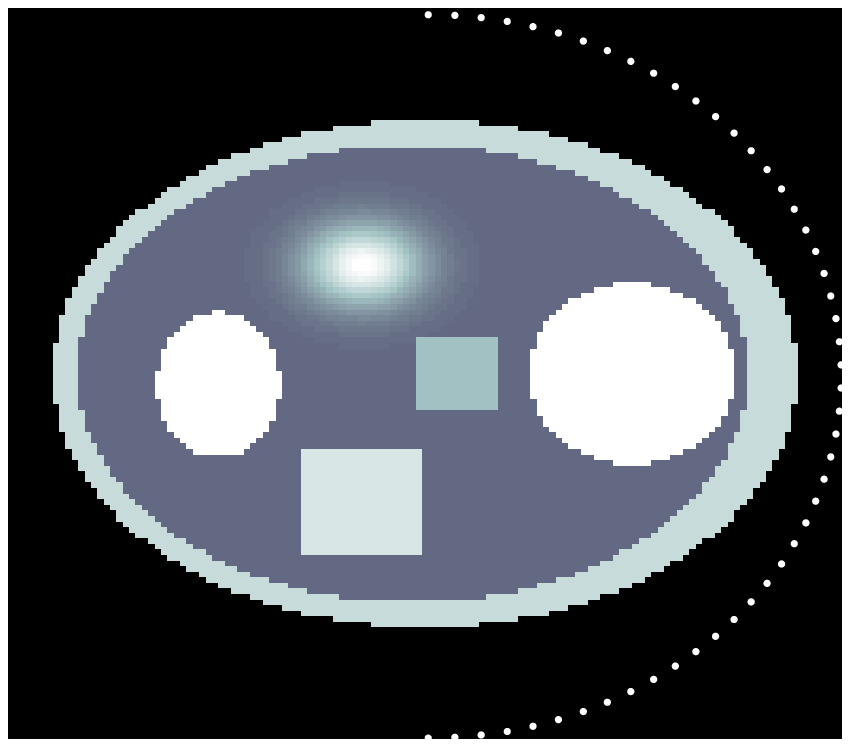}
\quad
\includegraphics[height = 0.4\textwidth, width=0.45\textwidth]{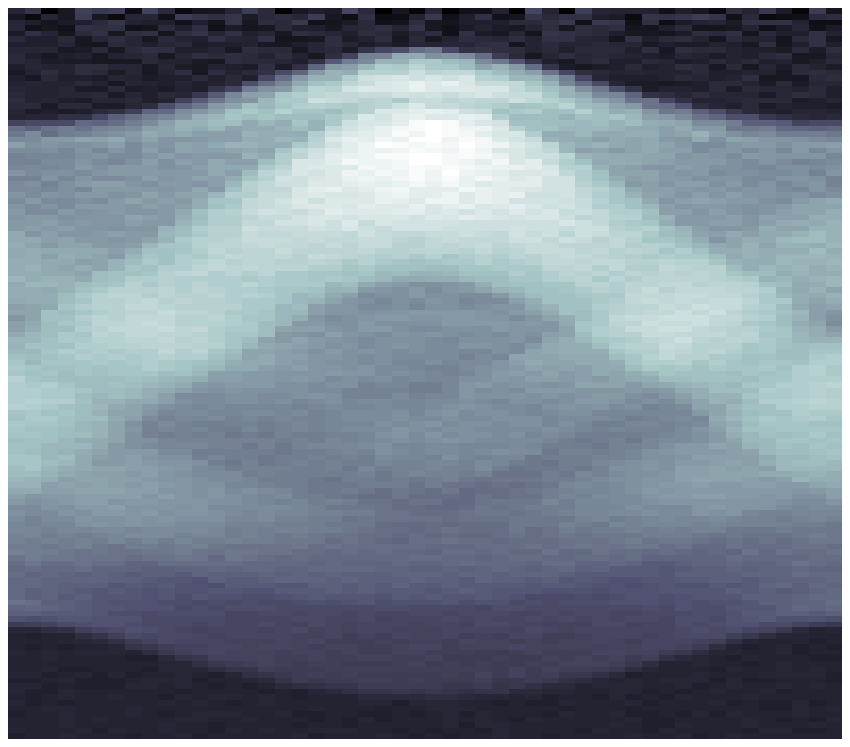}
\caption{The left picture shows the phantom $x^\dag$, where the white
dots indicate the locations of the detectors. The corresponding data
$(y_i^\delta)_i$ are depicted on the right.}
\label{fg:tat-original}
\end{figure}

\begin{psfrags}
\psfrag{LLK}{\color{white}\textsc{l-LK}}
\psfrag{LSDK}{\color{white}\textsc{l-SDK}}
\begin{figure}[b!]
\begin{center}
\includegraphics[width = 0.48\textwidth, height = 0.4 \textwidth]{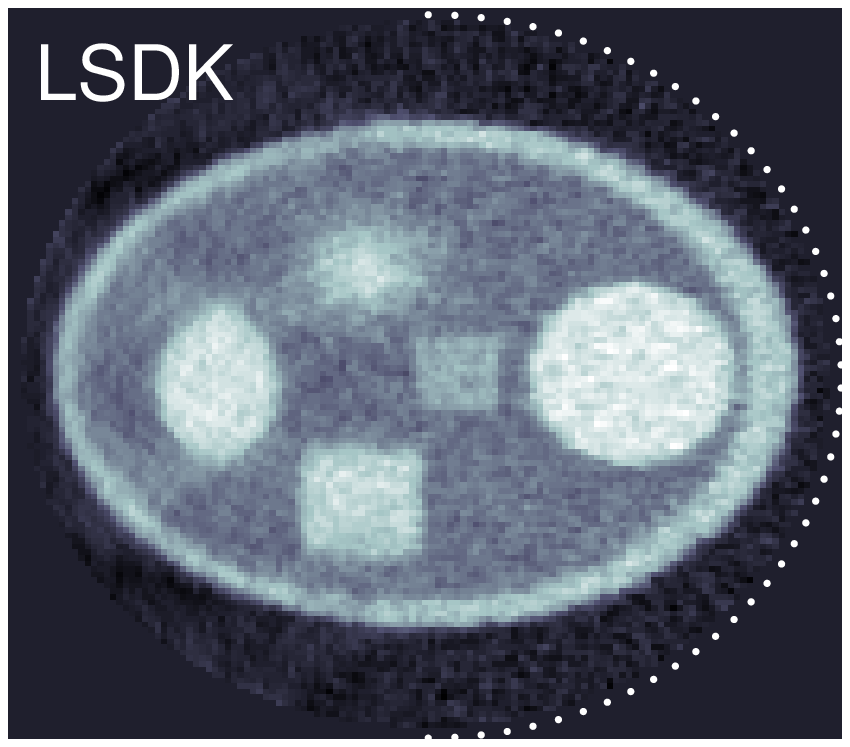}
\includegraphics[width = 0.48\textwidth, height = 0.4 \textwidth]{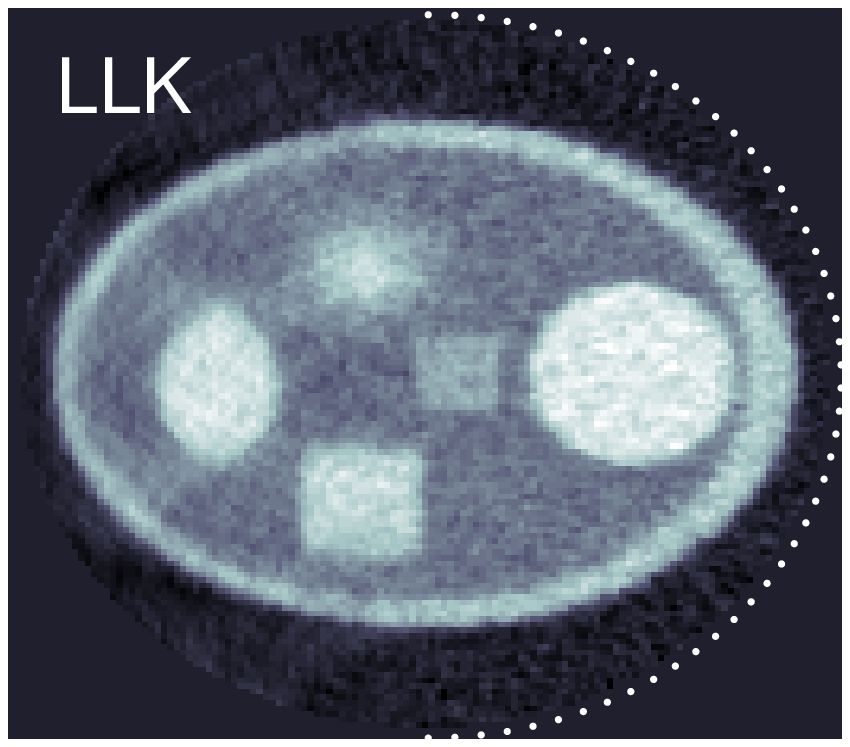}
\caption{Numerical reconstructions $x^\delta_{k^\delta_*}$ with
$\alpha_{\rm min} =0.4$ of the phantom depicted in Figure
\ref{fg:tat-original}.} \label{fg:ta-recon}
\end{center}
\end{figure}
\end{psfrags}

\medskip
In the following numerical examples, we consider the \textsc{l-SDK}
method with either the choice $\pr(s) = \min( \alpha_{\rm min} s ,
2) $ or $\pr(s) = \alpha_{\rm min}$ (which corresponds to the
\textsc{l-LK} method). In both cases we use $\alpha_{\rm min} = 0.4$
or $\alpha_{\rm min} = 1$, and   assume $N = 50$ measurements. The
phantom  $x^\dag$,  shown in the left picture in Figure
\ref{fg:tat-original},  consists  of a superposition of
characteristic functions and one Gaussian kernel.  Data $y_i = \Mo_i
x^\dag$ were calculated via numerical integration with the
trapezoidal rule and $4 \%$ noise was added, such that $\norm{y_i -
y_i^\delta}/\norm{y_i} \approx 0.04$.  The  the regularized
solutions $x^\delta_{k_*^\delta}$ with $\alpha_{\rm min} =0.4$ are
depicted in Figure \ref{fg:ta-recon}. For both, the \textsc{l-SDK}
\textsc{l-LK} method, all visible parts of the phantom $x^\dag$ are
reconstructed reliable.

\begin{psfrags}
\psfrag{w1}{\footnotesize $\alpha_{\rm min} = 0.4$}
\psfrag{w2}{\footnotesize $\alpha_{\rm min} = 1$} \psfrag{1}{ }
\psfrag{0.4}{ }
\begin{figure}[hbt!]
\begin{center}
\includegraphics[height = 0.3\textwidth]{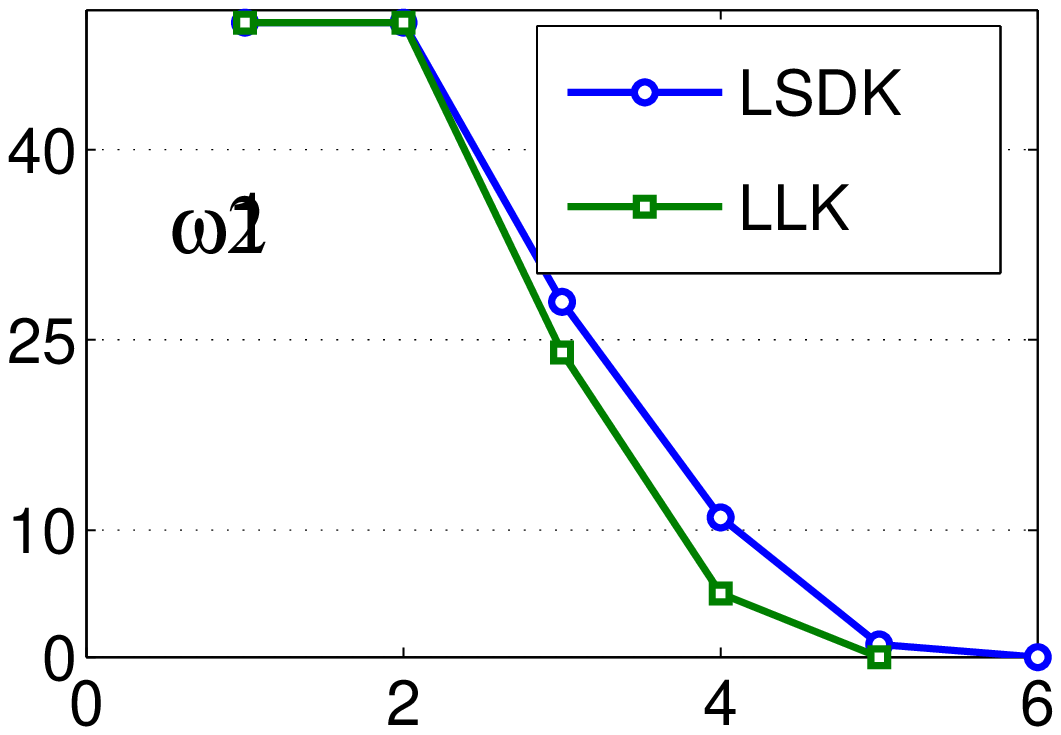}
\includegraphics[height = 0.3\textwidth]{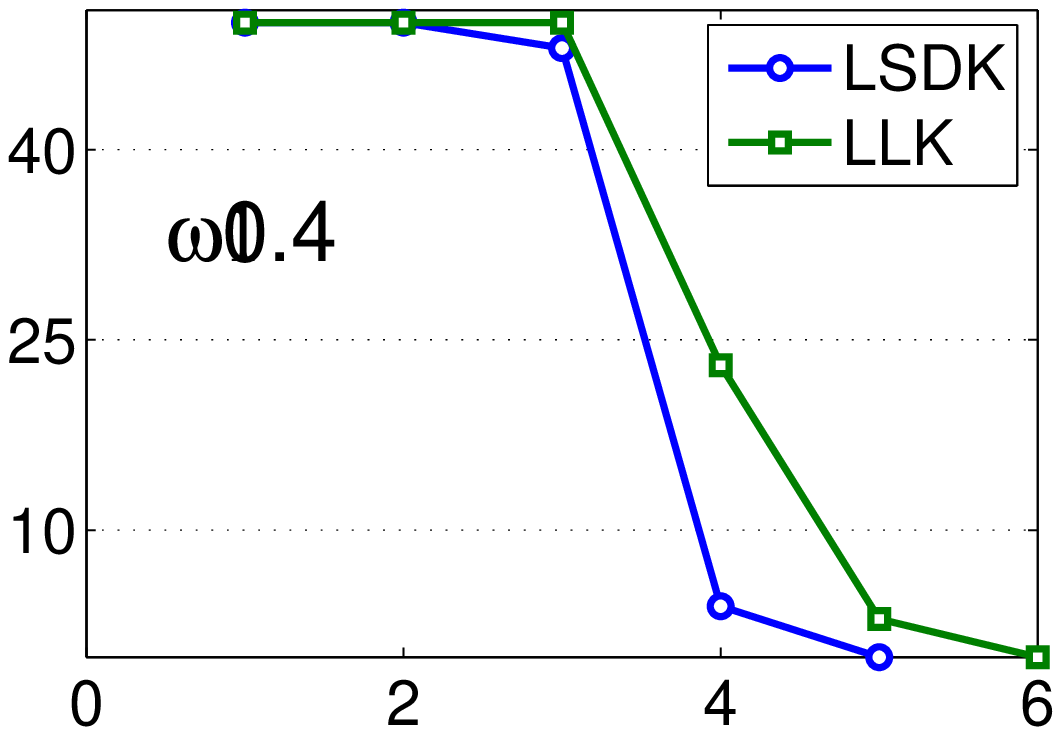}
\caption{The $x$-axis shows the number of cycles, while the number of
actually performed iterations within each cycle is shown at the $y$-axis.}
\label{fg:ta-steps}
\end{center}
\end{figure}
\end{psfrags}

\begin{psfrags}
\psfrag{w1}{\footnotesize $\alpha_{\rm min} = 0.4$}
\psfrag{w2}{\footnotesize $\alpha_{\rm min} = 1$} \psfrag{1}{ }
\psfrag{0.4}{ }
\begin{figure}[htb!]
\begin{center}
\includegraphics[height = 0.31 \textwidth]{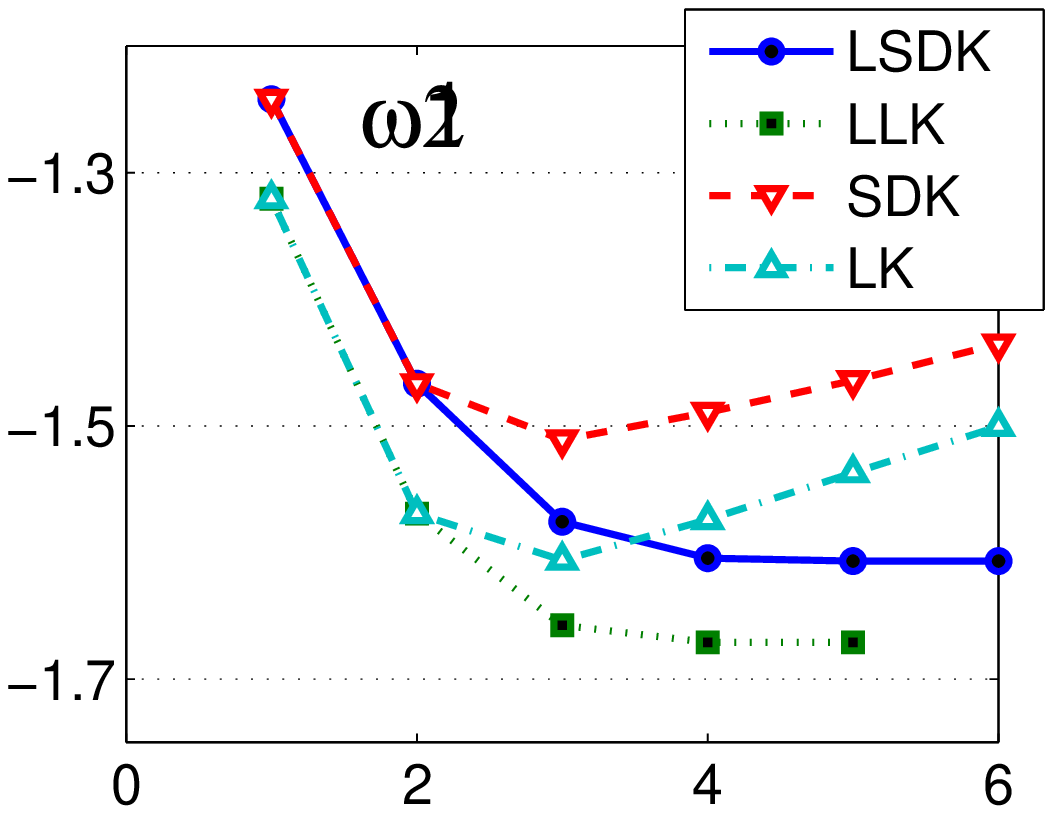}
\includegraphics[height = 0.31 \textwidth]{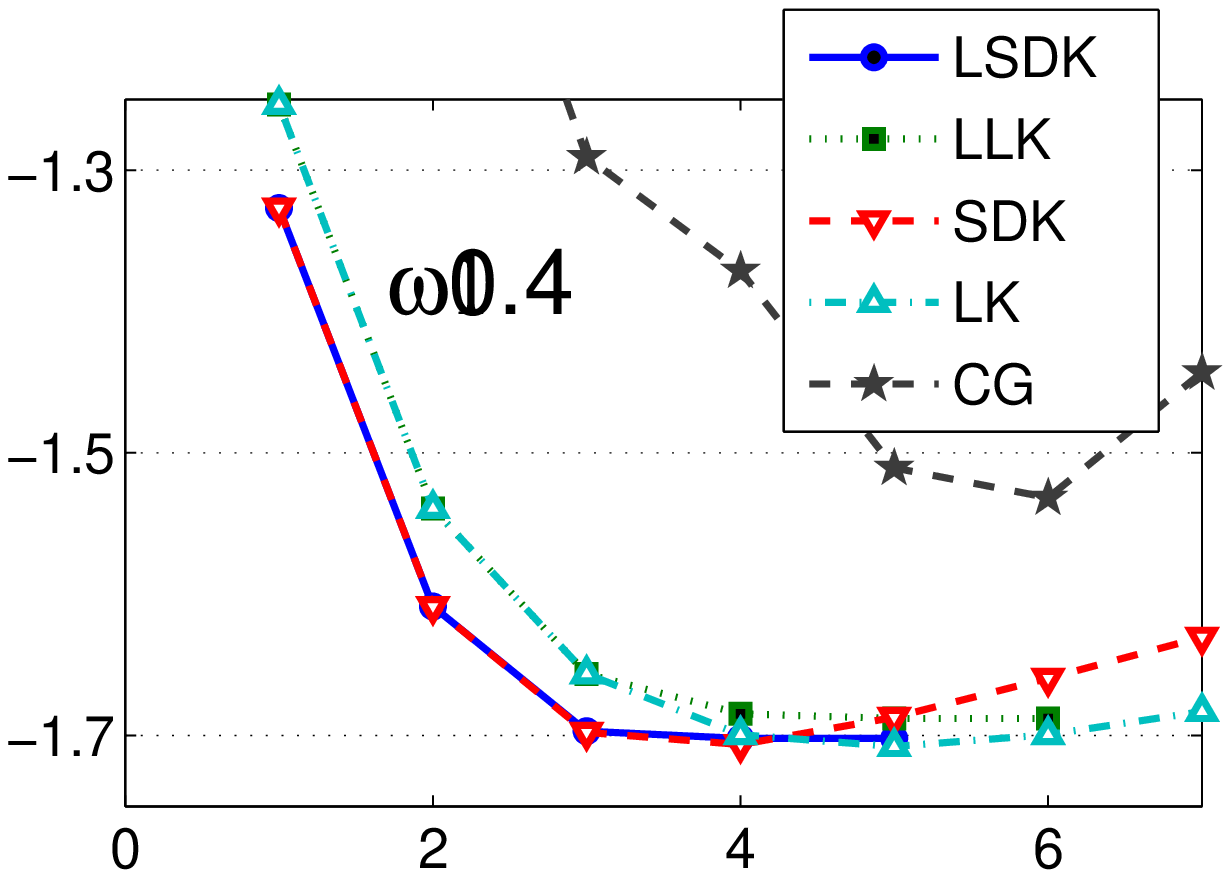}
\caption{ Evolution of the relative error $\ln
\norm{x^\dag-x_k^\delta}/\norm{x^\dag}$.}
\label{fg:ta-error}\end{center}
\end{figure}
\end{psfrags}

Figure \ref{fg:ta-steps} and Figure  \ref{fg:ta-error} show the
number of actually performed iterations and the reconstruction error
$e_k^\delta := \norm{x_k^\delta - x^\dag}$ respectively. For
comparison purposes, the error for the \textsc{SDK} and the
\textsc{LK} iteration (without loping parameter) are also included.
In all cases, the smaller relaxation parameter $\alpha_{\rm min}$
gives the smaller reconstruction errors.  This behavior is typically
for the application of Kaczmarz type iterations to Radon transforms
\cite{CenEggGor83, Nat01};  therefore in praxis often relatively
small relaxation parameters are chosen. For $\alpha_{\rm min} =1$,
the loping strategy significantly reduces the reconstruction error
of the no-loping iterations. Also, for $\alpha_{\rm min} =0.4$, the
regularized solution of the loping Kaczmarz methods (automatically
stopped according to \eqref{eq:def-discrep}) have errors comparable
to the optimal solution of their non-loping counterparts when
stopped after the cycle with minimal error (which is not available
in practice).

\begin{table}[b!]
\centering
\begin{tabular}{l | @{\hspace{1cm}} c @{\hspace{0.5cm}} c @{\hspace{0.5cm}} l }
\toprule
                  & Cycles & Runtime (sec) & Error ($\%$) \\
\midrule
\textsc{L-SDK}    & 5      & 21.9          & 18.2  \\
\textsc{L-LK}     & 6      & 21.4          & 18.5  \\
\textsc{SDK}      & 4      & 24.5          & 18.2  \\
\textsc{LK}       & 5      & 16.9          & 18.1  \\
\textsc{CGNE}     & 5      & 38.2          & 21.6  \\
\bottomrule
\end{tabular}
\caption{Comparison of the performance of different iterative methods.
The non-loping iterations are stopped after the cycle with minimal error.}
\label{tb:ta-summary}
\end{table}

To point out  the effectiveness of the  loping Kaczmarz methods for
solving linear inconsistent systems we included the reconstruction
error for the CGNE iteration (conjugate gradient \cite{HesSte53,
She94} applied to normal equations). If  stooped appropriately  the
CGNE method is known to be a regularization method \cite{EngHanNeu96, Han95}.
As can be seen in Figure \ref{fg:ta-error} the reconstruction error for
the \textsc{l-SDK} and the \textsc{l-LK} methods is much smaller that that
for the CGNE iteration.
In Table \ref{tb:ta-summary} run times for reconstructing an image on a
$120 \times 120$ grid  are compared (with non-optimized Matlab implementation
on  iMac with 2 GHz Intel Core Duo processor).


\section{An Inverse Doping Problem}
\label{sec:app-2}

In this section we present another comparison of the numerical
performance of the \textsc{l-SDK}, \textsc{l-LK} and \textsc{LK}
methods. This time we consider an application related to {\em
inverse doping problems} for semiconductors \cite{BELM04, HKLS07,
LMZ06, BurEtAl01} For details on the mathematical modeling of this
inverse problem we refer the reader to \cite[Section~3]{HKLS07}.
\medskip

In what follows we describe the abstract formulation in Hilbert
spaces of the problem (the so called {\em inverse doping problem in
the linearized unipolar model for current flow measurements}). Let
$\Omega := (0,1) \times (0,1) \subset \mathbb{R}^2$ be the domain
representing the semiconductor device (a diode). The two
semiconductor contacts are represented by the boundary parts:
$$
\Gamma_0  :=  \set{ (s,0) :   \ s \in (0,1)  } \, , \quad \Gamma_1
:=  \set{ (s,1) : \ s \in (0,1)    } \, ,
$$
(we denote $\partial\Omega_D := \Gamma_0 \cup \Gamma_1$) while the
insulated surfaces of the semiconductor are represented by
$\partial\Omega_N  :=  \set{ (0,t) : \ t \in (0,1) } \cup \set{
(1,t) : \ t \in (0,1) }$. This specific inverse doping problem can
be reduced to the identification of the positive parameter function
$x$ (the doping profile $C$ is related to $x$ by $C = x - \lambda^2
\Delta (\ln x)$) in the model
\begin{eqnarray} \label{eq:upolU-A}
    \mu_n \, \nabla \cdot ( x(\xi) \, \nabla u)
    = 0 \, , \quad           & \mbox{in } \Omega \\ \label{eq:upolU-B}
    u
    = U(\xi) \, , \quad        & \mbox{on } \partial\Omega_D \\ \label{eq:upolU-C}
    \nabla u \cdot\nu
    = 0 \, , \quad           & \mbox{on } \partial\Omega_N
\end{eqnarray}
from measurements of the Voltage--Current map (the forward operator)
\begin{eqnarray*}
    \Sigma_x: H^{3/2}(\partial\Omega_D)  & \to &  \mathbb R \, , \\
    U  & \mapsto &  \mu_n \int_{\Gamma_1} e^{V_{\rm bi}(\xi)} u_\nu(\xi) \, d\Gamma
\end{eqnarray*}
which maps an applied potential $U$ at $\partial\Omega_D$ to the
corresponding \textit{total current flow} $\Sigma_x(U)$ through the
contact $\Gamma_1$. Here $\mu_n$, $\lambda$ are positive constants
and $V_{\rm bi}$ is a known logarithmic function defined on
$\partial\Omega_D$.

\begin{figure}[tb!]
\centerline{
\includegraphics[width=0.50\textwidth]{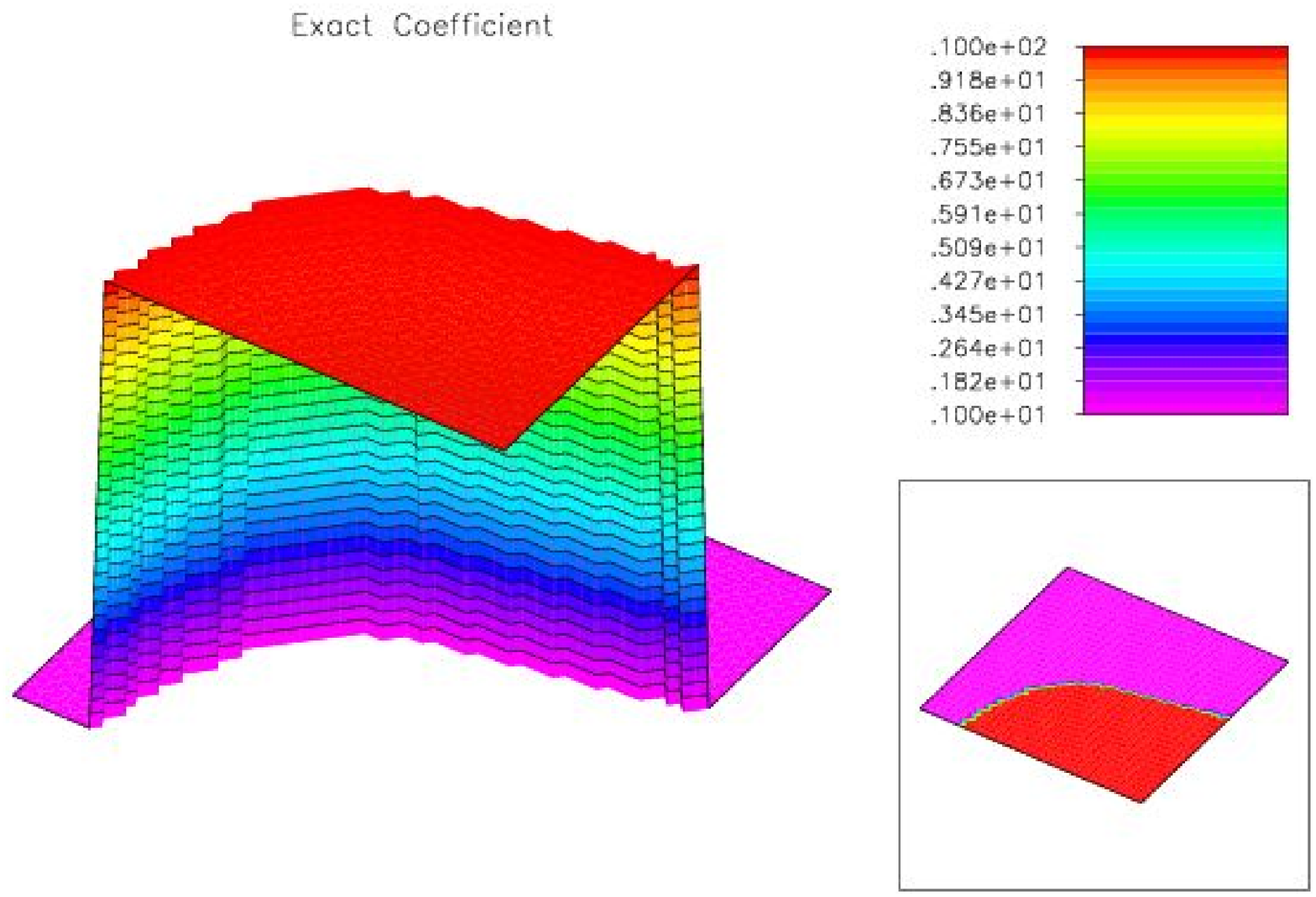} \hskip0.4cm
\includegraphics[width=0.50\textwidth]{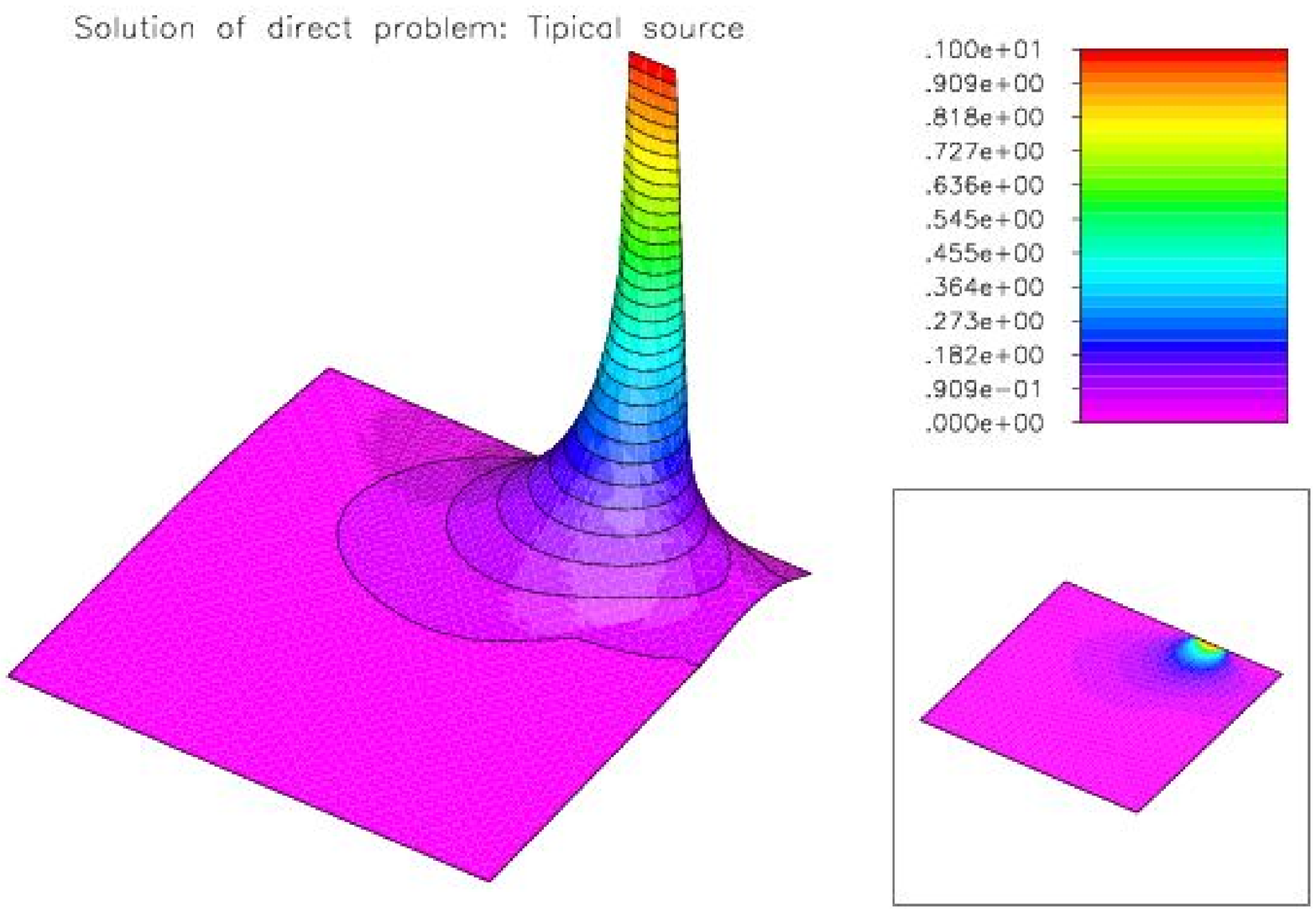} }
\centerline{
\includegraphics[width=0.50\textwidth]{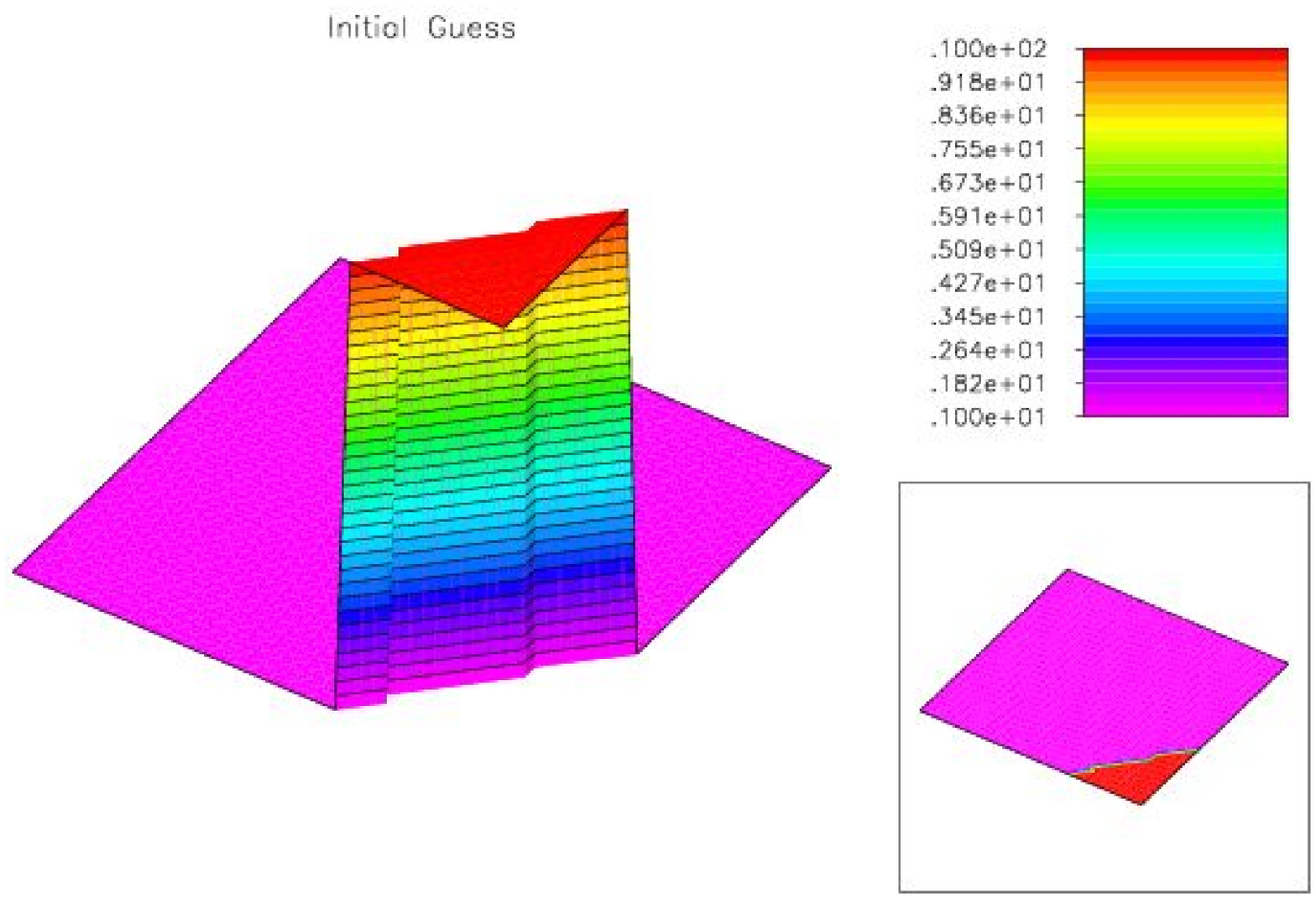}}
\caption{\small In the top left picture, the doping profile to be
identified. In the top right picture, a typical voltage profile
$U_i$ and the corresponding solution $u$ of
\eqref{eq:upolU-A} - \eqref{eq:upolU-C}. The initial guess used for
the \textsc{l-SDK}, \textsc{l-LK} and \textsc{LK} iterative methods
is shown in the bottom picture. The boundary parts $\Gamma_0$ and
$\Gamma_1$ correspond to the top right and to the lower left edge
respectively (the origin is the right corner).}
\label{fig:idop-setup}
\end{figure}

Due to the nature of the practical experiments that can be performed
on a factory environment, some restrictions on the data have to be
taken into account:
\begin{enumerate}
\item \label{item1}
The voltage profiles $U \in H^{3/2}(\partial\Omega_D)$ must satisfy
$U(\xi) = 0$ at the contact $\Gamma_1$.
\item The parameter $x$ has to be determined from a finite number of
measurements, i.e. from the data
\begin{equation} \label{eq:data-semicond}
    y_i^\delta := \Sigma_x(U_i) \in Y:= \mathbb R\, ,
    \quad i = 0, \dots, N-1 \, ,
\end{equation}
where the $U_i \in H^{3/2}(\partial\Omega_D)$ are prescribed voltage
profiles satisfying Item \ref{item1}.
\end{enumerate}
Therefore we can model the inverse doping problem with a system of
operator equations of the form \eqref{eq:inv-probl}, namely
$$
     \F_i(x) = y_i^\delta \, , \quad i = 0, \dots , N-1 \, ,
$$
where $x \in L^2(\Omega) =: X$ is the unknown parameter, $y_i^\delta
\in \R =: Y$ denote the measured data, $\F_i: X \to Y$ defined by
$\F_i(x) := \Sigma_x(U_i)$ are the parameter to output maps, with
domains of definition
$$
    D_i := \{ x \in L^\infty(\Omega) :  \, 0 < x_{\rm min}
   \le x  \le x_{\rm max}, \mbox{ a.e.} \} \, .
$$
It is worth mentioning that, although the operators $\F_i$ are
Fr\'echet differentiable, they do not satisfy the tangential cone
condition \eqref{eq:a-tcc}. Therefore, the convergence results
derived in Section~\ref{sec:conv-an} cannot be applied.

\begin{figure}[tb!]
\centerline{
\includegraphics[width=0.49\textwidth]{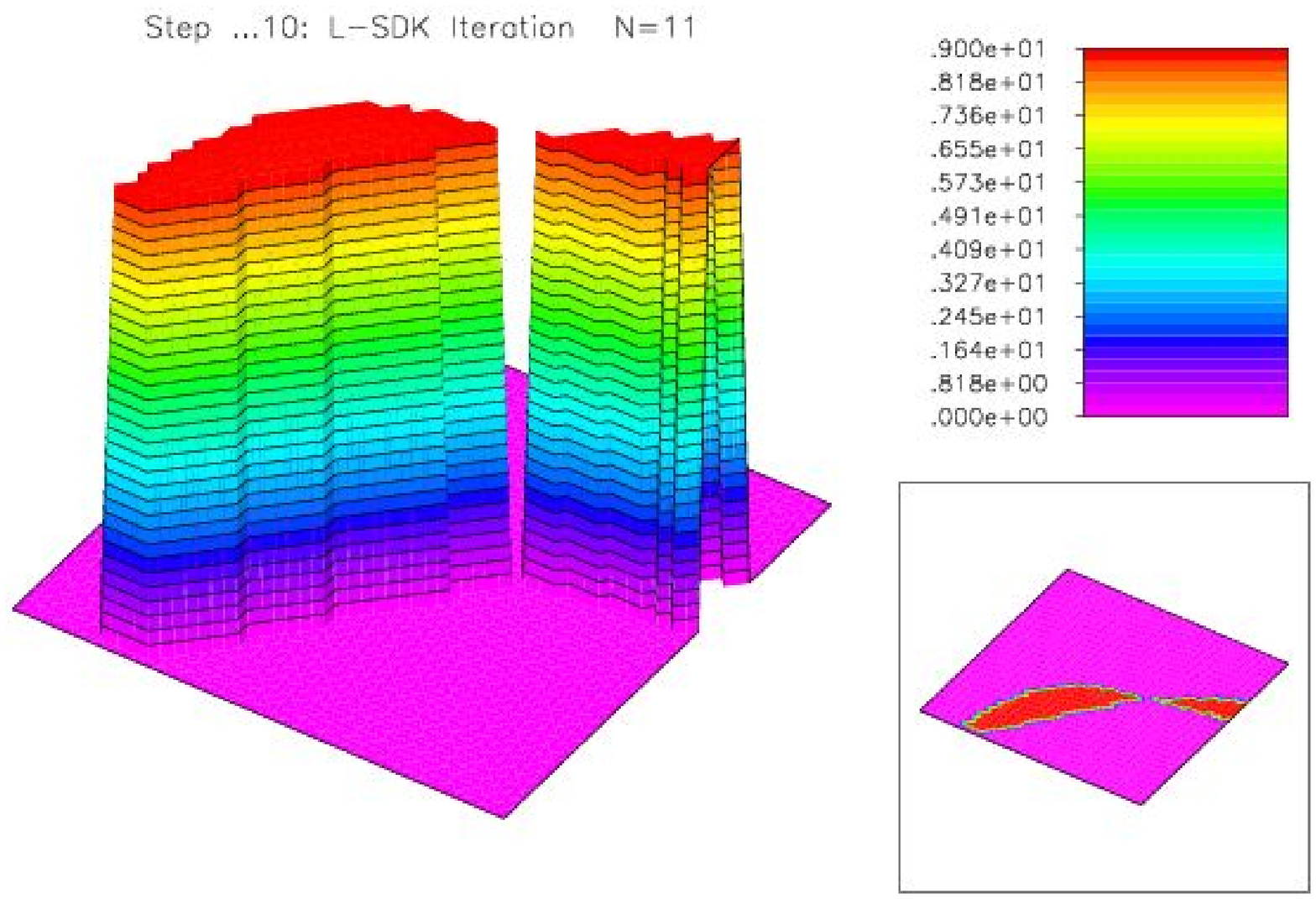} \hskip0.4cm
\includegraphics[width=0.49\textwidth]{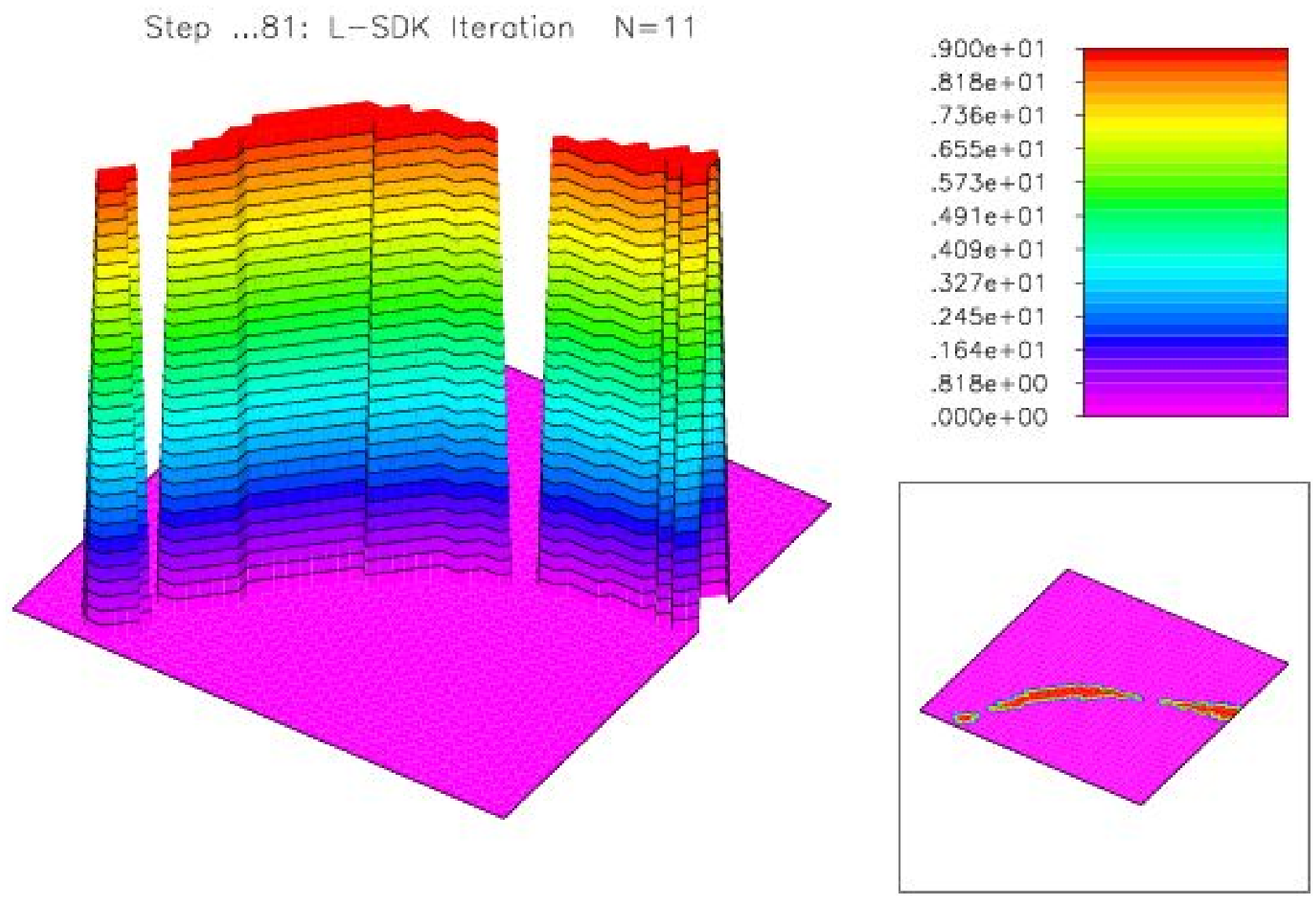} }
\bigskip
\centerline{
\includegraphics[width=0.49\textwidth]{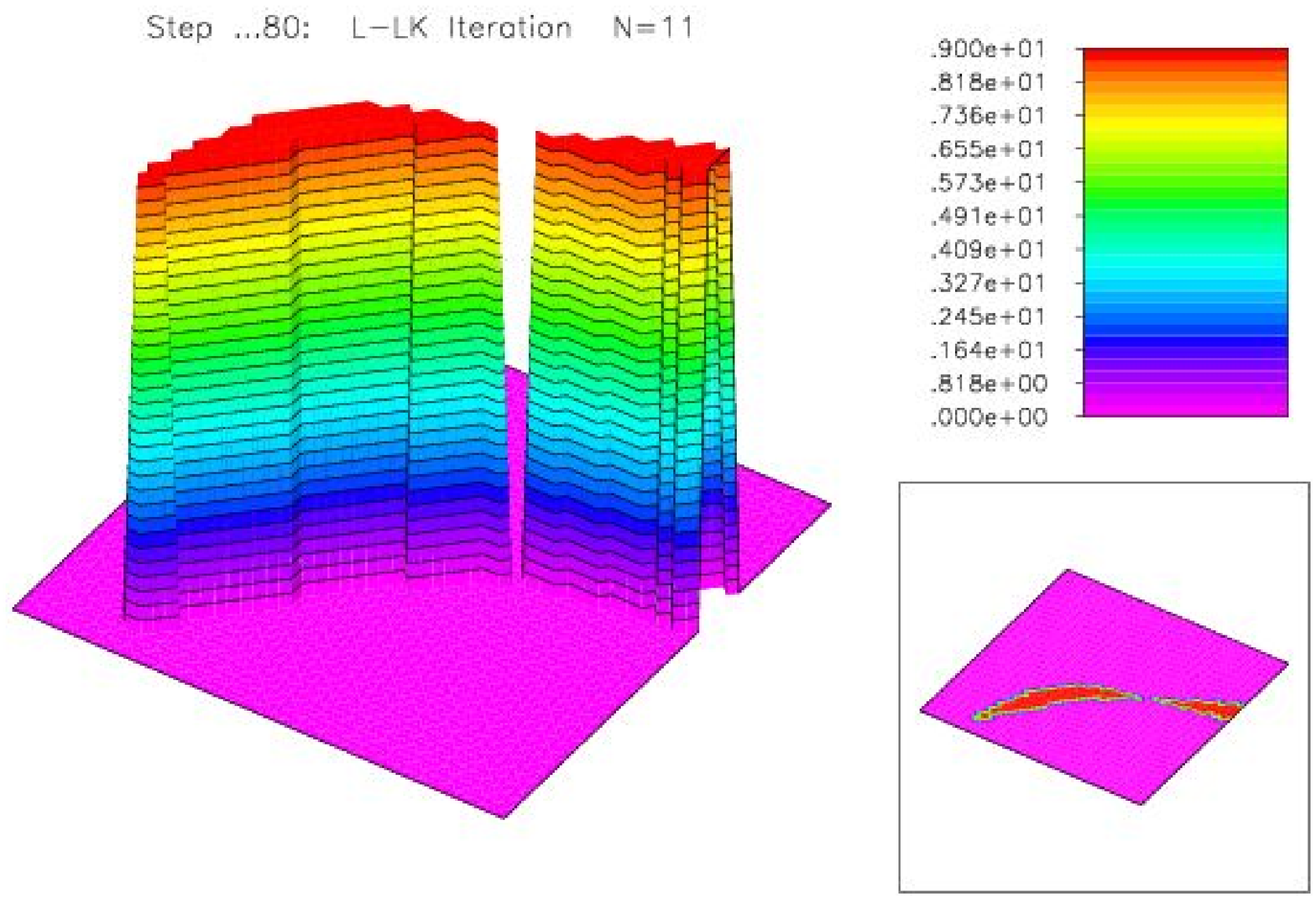} \hskip0.4cm
\includegraphics[width=0.49\textwidth]{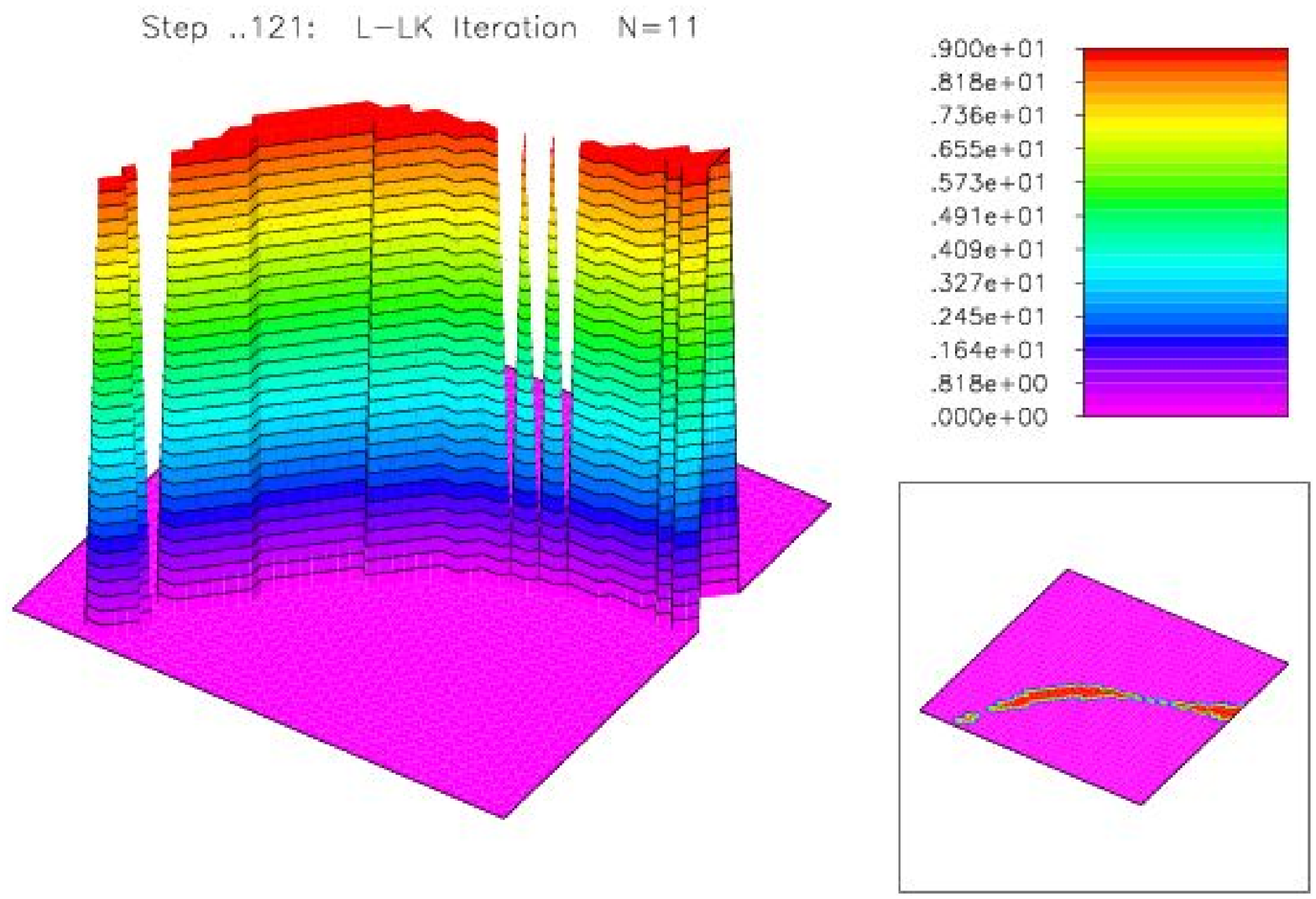} }
\bigskip
\centerline{
\includegraphics[width=0.49\textwidth]{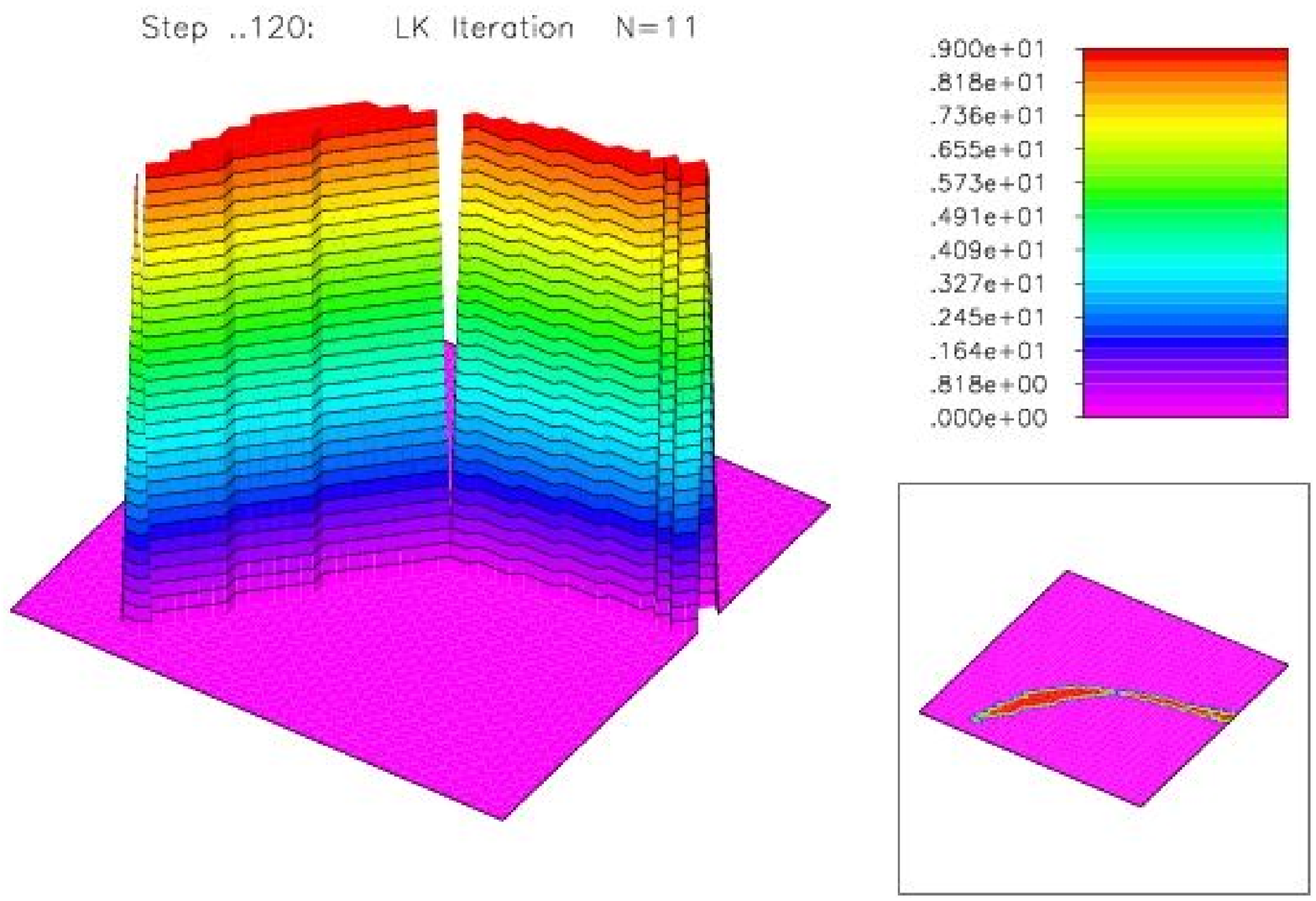} \hskip0.4cm
\includegraphics[width=0.49\textwidth]{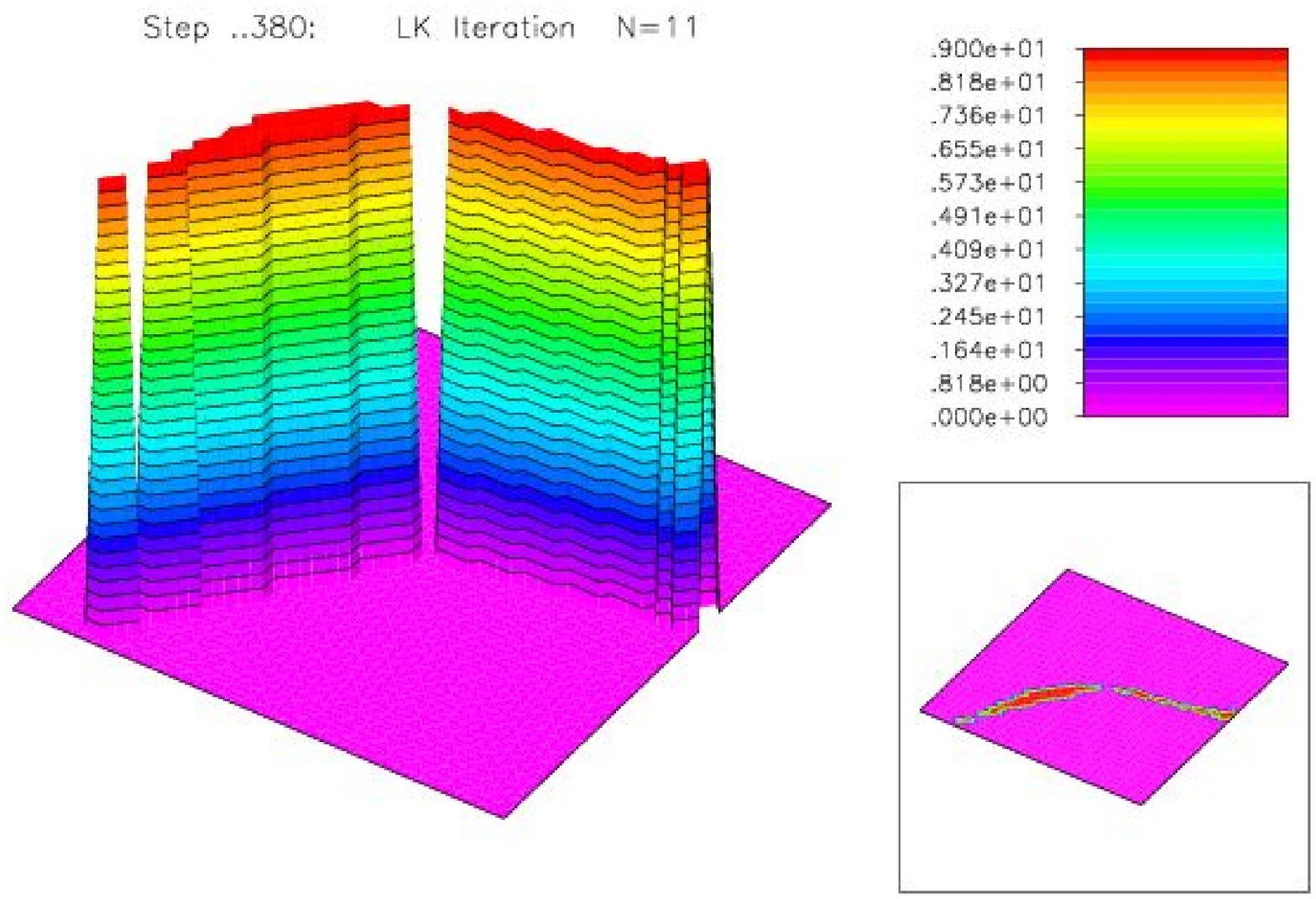} }
\caption{\small Comparison between the \textsc{l-SDK}, \textsc{l-LK}
and \textsc{LK} methods. The top two pictures show the iterative errors
obtained by the \textsc{l-SDK} iteration after 10 and 81 cycles. The
two central pictures show the iterative errors obtained by the \textsc{l-LK}
iteration after 80 and 121 cycles. The two bottom pictures show the iterative
errors obtained by the \textsc{LK} after 120 and 380 cycles.}
\label{fig:idop-iter}
\end{figure}

\medskip
In the following numerical examples we assume that $N = 11$
Dirichlet--Neumann pairs $(U_i, \F_i(x') )$ of measurement data are
available. The fixed inputs $U_i$, are chosen to be piecewise
constant functions supported in $\Gamma_0$,
$$
  U_i(s) \ := \ \left\{ \begin{array}{rl}
      1 & |s - s_i| \le h \\
      0 & {\rm else} \end{array} \right. \,,
$$
where the points $s_i$ are uniformly distributed on $\Gamma_0$ and
$h = 1/32$. The doping profile to be reconstructed is shown in
Figure~\ref{fig:idop-setup} (top left picture). The top right
picture of Figure~\ref{fig:idop-setup} shows a typical voltage
profile $U_j$ (applied at $\Gamma_0$) as well as the corresponding
solution $u$ of \eqref{eq:upolU-A} - \eqref{eq:upolU-C}. In these
pictures, as well as in the forthcoming ones, $\Gamma_1$ appears on
the lower left edge and $\Gamma_0$ on the top right edge (the origin
corresponds to the upper right corner).

In Figure~\ref{fig:idop-iter} we show the evolution of the iteration
error for the \textsc{l-SDK}, \textsc{l-LK} and \textsc{LK} methods.
The same initial guess was used for the three methods (see
Figure~\ref{fig:idop-setup}). In our computations we chose $\tau =
2.5$ in \eqref{eq:def-tau}. The stopping rule for the \textsc{l-SDK}
method is satisfied after 81 cycles. For the \textsc{l-LK} method,
the same stopping criteria is reached only after 121 cycles. In
order to obtain the same accuracy with the \textsc{LK} method, 380
cycles are required. In the top pictures of
Figure~\ref{fig:idop-iter} one can see the iteration error for the
\textsc{l-SDK} method after 10 and 81 cycles. For comparison
purposes, the iteration error for the \textsc{l-LK} method is shown
after 80 and 121 cycles (see the central pictures of
Figure~\ref{fig:idop-iter}). The bottom pictures of
Figure~\ref{fig:idop-iter} show the iteration error for the
\textsc{LK} method after 120 and 380 cycles. The number of actually
computed iterative steps within each cycle of the \textsc{l-SDK} and
\textsc{l-LK} methods is shown in Figure~\ref{fig:idop-loping}.

\begin{psfrags}
\psfrag{Computed iterative steps within each
cycle}{{}\hspace{0.1\textwidth}\small Computed iterative steps
within each cycle} \psfrag{20}{\footnotesize 20}
\psfrag{40}{\footnotesize 40} \psfrag{80}{\footnotesize 80}
\psfrag{60}{\footnotesize 60} \psfrag{100}{\footnotesize 100}
\psfrag{120}{\footnotesize 120} \psfrag{1}{\footnotesize 1}
\psfrag{3}{\footnotesize 3} \psfrag{5}{\footnotesize 5}
\psfrag{9}{\footnotesize 9} \psfrag{7}{\footnotesize 7}
\psfrag{11}{\footnotesize 11}
\begin{figure}[b!]
\centerline{
\includegraphics[height=0.45\textwidth,width=0.65\textwidth]{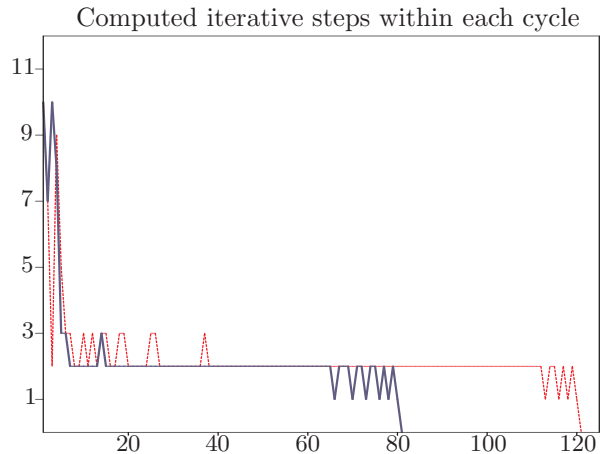} }
\caption{\small Comparison between the performance of \textsc{l-SDK}
and \textsc{l-LK} methods. The solid (blue) line shows the actually
performed number of steps within each cycle of the \textsc{l-SDK}
method, while the dashed (red) line gives the corresponding
information with respect to the \textsc{l-LK} method.}
\label{fig:idop-loping}
\end{figure}
\end{psfrags}

As one can see in Figure~\ref{fig:idop-loping}, no more than 2
steepest descent steps per cycle are computed after the $14$-th
cycle of the \textsc{l-SDK} method. Analogously, no more than 2
Landweber steps per cycle are computed after the $37$-th cycle of
the \textsc{l-LK} method. In total, for the computation of the
\textsc{LK}-approximation in Figure~\ref{fig:idop-iter} (380
cycles), 4180 Landweber steps are needed, while the
\textsc{l-LK}-approximation (121 cycles) requires the computation of
258 Landweber steps and the \textsc{l-SDK}-approximation (81 cycles)
requires the computation of 184 steepest descent steps. The
\textsc{l-LK} method requires almost 50\% more cycles than the
\textsc{l-SDK} method in order to reach the stopping criteria
\eqref{eq:def-discrep}. Moreover, the \textsc{LK} method requires
almost three times more cycles than the \textsc{l-LK} method in
order to achieve the same accuracy (see \cite{HKLS07} for other
comparisons between the \textsc{LK} and \textsc{l-LK} methods).

The efficiency of the \textsc{l-SDK} method becomes even more
evident when we compare the total number of actually performed
iterative steps. Each cycle of the \textsc{LK} method requires the
computation of 11 steps, while in the \textsc{l-SDK} and
\textsc{l-LK} methods the number of actually performed steps per
cycle is very small after a few number of cycles.

\section{Conclusions} \label{sec:conclusion}

In this paper we propose a new iterative method for inverse problems
of the form \eqref{eq:inv-probl}, namely the \textsc{l-SDK} method.
As a by-product we also formulated the \textsc{SDK} iteration, which
is the steepest descent counterpart of the \textsc{LK} method
\cite{KowSch02}. In the \textsc{l-SDK} iteration we omit an update
of the \textsc{SDK} iteration (within one cycle) if corresponding
$i$-th residual is below some threshold. Consequently, the
\textsc{l-SDK} method is not stopped until all residuals are below
the specified threshold. We provided a complete convergence analysis
for the \textsc{l-SDK} iteration, proving that it is a convergent
regularization method in the sense of \cite{EngHanNeu96}.

The abstract theory was applied to thermoacoustic computed
tomography and an inverse problem for semiconductors. In both
applications the \textsc{l-SDK} method turned out to be an efficient
iterative regularization method.

\section*{Acknowledgments}

The work of M.H. and O.S. is supported by the FWF (Austrian Science
Fund) grants Y--123INF and P18172--N02. The work of A.L. and A.DC.
are supported by the Brazilian National Research Council CNPq,
grants 306020/2006--8 and 474593/2007--0. The authors thank Andreas
Rieder for stimulating discussion on iterative regularization
methods.
\bibliographystyle{amsplain}
\bibliography{steepd-kac}

\end{document}